\input amstex
\documentstyle{amsppt}
\input epsf.tex
%\input label.def
%\input debug.def

%Version: Alex Degtyarev (home) 28/09/98  16:11:30
\def\stydate{May 10, 2002}

\chardef\tempcat\catcode`\@ \ifx\undefined\amstexloaded\input
amstex \else\catcode`\@\tempcat\fi \expandafter\ifx\csname
amsppt.sty\endcsname\relax\input amsppt.sty \fi
\let\tempcat\undefined

\immediate\write16{This is LABEL.DEF by A.Degtyarev <\stydate>}
\expandafter\ifx\csname label.def\endcsname\relax\else
  \message{[already loaded]}\endinput\fi
\expandafter\edef\csname label.def\endcsname{%
  \catcode`\noexpand\@\the\catcode`\@\edef\noexpand\styname{LABEL.DEF}%
  \def\expandafter\noexpand\csname label.def\endcsname{\stydate}%
    \toks0{}\toks2{}}
\catcode`\@11
\def\labelmesg@ {LABEL.DEF: }
{\edef\temp{\the\everyjob\W@{\labelmesg@<\stydate>}}
\global\everyjob\expandafter{\temp}}

%%%%%%%%%%%%%%%%%%%% Handy stuff
\def\@car#1#2\@nil{#1}
\def\@cdr#1#2\@nil{#2}
\def\eat@bs{\expandafter\eat@\string}
\def\eat@ii#1#2{}
\def\eat@iii#1#2#3{}
\def\eat@iv#1#2#3#4{}
\def\@DO#1#2\@{\expandafter#1\csname\eat@bs#2\endcsname}
\def\@N#1\@{\csname\eat@bs#1\endcsname}
\def\@Nx{\@DO\noexpand}
\def\@Name#1\@{\if\@undefined#1\@\else\@N#1\@\fi}
\def\@Ndef{\@DO\def}
\def\@Ngdef{\global\@Ndef}
\def\@Nedef{\@DO\edef}
\def\@Nxdef{\global\@Nedef}
\def\@Nlet{\@DO\let}
\def\@undefined#1\@{\@DO\ifx#1\@\relax\true@\else\false@\fi}
\def\@@addto#1#2{{\toks@\expandafter{#1#2}\xdef#1{\the\toks@}}}
\def\@@addparm#1#2{{\toks@\expandafter{#1{##1}#2}%
    \edef#1{\gdef\noexpand#1####1{\the\toks@}}#1}}
\def\make@letter{\edef\t@mpcat{\catcode`\@\the\catcode`\@}\catcode`\@11 }
\def\donext@{\expandafter\egroup\next@}
\def\x@notempty#1{\expandafter\notempty\expandafter{#1}}
%%%%% Defines #1 to be #2 in lower case
\def\lc@def#1#2{\edef#1{#2}%
    \lowercase\expandafter{\expandafter\edef\expandafter#1\expandafter{#1}}}
%%%%% Finds a text in a comma separated list
\newif\iffound@
\def\find@#1\in#2{\found@false
    \DNii@{\ifx\next\@nil\let\next\eat@\else\let\next\nextiv@\fi\next}%
    \edef\nextiii@{#1}\def\nextiv@##1,{%
    \edef\next{##1}\ifx\nextiii@\next\found@true\fi\FN@\nextii@}%
    \expandafter\nextiv@#2,\@nil}
%%%%% Disable \outer'ness
{\let\head\relax\let\specialhead\relax\let\subhead\relax
\let\subsubhead\relax\let\proclaim\relax
\gdef\let@relax{\let\head\relax\let\specialhead\relax\let\subhead\relax
    \let\subsubhead\relax\let\proclaim\relax}}
%%%%% Hacks
\newskip\@savsk
% We add a tiny skip in order to make \@esphack work correct
\let\@ignorespaces\ignorespaces
\def\@ignorespacesp{\ifhmode
  \ifdim\lastskip>\z@\else\penalty\@M\hskip-1sp%
        \penalty\@M\hskip1sp \fi\fi\@ignorespaces}
\def\ignorespaces{\protect\@ignorespacesp}
% However, we have to redefine some control sequences
%%??\def~{\unskip\nobreak\ \@ignorespaces}
\def\@bsphack{\relax\ifmmode\else\@savsk\lastskip
  \ifhmode\edef\@sf{\spacefactor\the\spacefactor}\fi\fi}
\def\@esphack{\relax
  \ifx\penalty@\penalty\else\penalty\@M\fi   % if this is after \nobreak
  \ifmmode\else\ifhmode\@sf{}\ifdim\@savsk>\z@\@ignorespacesp\fi\fi\fi}
%%%%% Check for a *, \nofrills, and []
\let\@frills@\identity@
\let\@txtopt@\identyty@
\newif\if@star
\newif\if@write\@writetrue
\def\@numopt@{\if@star\expandafter\eat@\fi}
\def\checkstar@#1{\DN@{\@writetrue
  \ifx\next*\DN@####1{\@startrue\checkstar@@{#1}}%
      \else\DN@{\@starfalse#1}\fi\next@}\FN@\next@}
\def\checkstar@@#1{\DN@{%
  \ifx\next*\DN@####1{\@writefalse#1}%
      \else\DN@{\@writetrue#1}\fi\next@}\FN@\next@}
\def\checkfrills@#1{\DN@{%
  \ifx\next\nofrills\DN@####1{#1}\def\@frills@####1{####1\nofrills}%
      \else\DN@{#1}\let\@frills@\identity@\fi\next@}\FN@\next@}
\def\checkbrack@#1{\DN@{%
    \ifx\next[\DN@[####1]{\def\@txtopt@########1{####1}#1}%
    \else\DN@{\let\@txtopt@\identity@#1}\fi\next@}\FN@\next@}
\def\check@therstyle#1#2{\bgroup\DN@{#1}\ifx\@txtopt@\identity@\else
        \DNii@##1\@therstyle{}\def\@therstyle{\DN@{#2}\nextii@}%
    \expandafter\expandafter\expandafter\nextii@\@txtopt@\@therstyle.\@therstyle
    \fi\donext@}

%%%%%%%%%%%%%%%%%%%% \@input, \include, .aux files, etc.
\newread\@inputcheck
\def\@input#1{\openin\@inputcheck #1 \ifeof\@inputcheck \W@
  {No file `#1'.}\else\closein\@inputcheck \relax\input #1 \fi}

\def\loadstyle#1{\edef\next{#1}%
    \DN@##1.##2\@nil{\if\notempty{##2}\else\def\next{##1.sty}\fi}%
    \expandafter\next@\next.\@nil\lc@def\next@\next
    \expandafter\ifx\csname\next@\endcsname\relax\input\next\fi}

\let\pagebody@\pagebody
\let\pagetop@\empty
\let\pagebot@\empty
\let\@Xend\empty
\def\pagebody{\pagetop@\pagebody@\pagebot@\@Xend}
\let\@Xclose\empty

\newwrite\@Xmain
\newwrite\@Xsub
\def\W@X{\write\@Xout}
\def\make@Xmain{\global\let\@Xout\@Xmain\global\let\end\endmain@
  \xdef\@Xname{\jobname}\xdef\@inputname{\jobname}}
\begingroup
\catcode`\(\the\catcode`\{\catcode`\{12
\catcode`\)\the\catcode`\}\catcode`\}12
\gdef\W@count#1((\lc@def\@tempa(#1)%
    \def\\##1(\W@X(\global##1\the##1))%
    \edef\@tempa(\W@X(%
        \string\expandafter\gdef\string\csname\space\@tempa\string\endcsname{)%
        \\\pageno\\\cnt@toc\\\cnt@idx\\\cnt@glo\\\footmarkcount@
        \@Xclose\W@X(}))\expandafter)\@tempa)
\endgroup
\def\readaux{\bgroup\checkbrack@\readaux@}
\let\begin\readaux
\def\readaux@{%
    \W@{>>> \labelmesg@ Run this file twice to get x-references right}%
    \global\everypar{}%
    {\def\\{\global\let}%
        \def\/##1##2{\gdef##1{\wrn@command##1##2}}%
        \disablepreambule@cs}%
    \make@Xmain{\make@letter\setboxz@h{\@input{\@txtopt@{\@Xname.aux}}%
            \lc@def\@tempa\jobname\@Name\open@\@tempa\@}}%
  \immediate\openout\@Xout\@Xname.aux%
    \immediate\W@X{\relax}\egroup}
\everypar{\global\everypar{}\readaux}
{\toks@\expandafter{\topmatter}
\global\edef\topmatter{\noexpand\readaux\the\toks@}}
\let\@@end@@\end

\def\@Xclose@{{\def\@Xend{\ifnum\insertpenalties=\z@
        \W@count{close@\@Xname}\closeout\@Xout\fi}%
    \vfill\supereject}}
\def\endmain@{\@Xclose@
    \W@{>>> \labelmesg@ Run this file twice to get x-references right}%
    \@@end@@}
\def\disablepreambule@cs{\\\disablepreambule@cs\relax}

\def\include#1{\bgroup
  \ifx\@Xout\@Xsub\DN@{\errmessage
        {\labelmesg@ Only one level of \string\include\space is supported}}%
    \else\edef\@tempb{#1}\clearpage
      \DN@##1 {\if\notempty{##1}\edef\@tempb{##1}\DN@####1\eat@ {}\fi\next@}%
    \DNii@##1.{\edef\@tempa{##1}\DN@####1\eat@.{}\next@}%
        \expandafter\next@\@tempb\eat@{} \eat@{} %
    \expandafter\nextii@\@tempb.\eat@.%
        \relaxnext@
      \if\x@notempty\@tempa
          \edef\nextii@{\write\@Xmain{%
            \noexpand\string\noexpand\@input{\@tempa.aux}}}\nextii@
        \ifx\undefined\@includelist\found@true\else
                    \find@\@tempa\in\@includelist\fi
            \iffound@\ifx\undefined\@noincllist\found@false\else
                    \find@\@tempb\in\@noincllist\fi\else\found@true\fi
            \iffound@\lc@def\@tempa\@tempa
                \if\@undefined\close@\@tempa\@\else\edef\next@{\@Nx\close@\@tempa\@}\fi
            \else\xdef\@Xname{\@tempa}\xdef\@inputname{\@tempb}%
                \W@count{open@\@Xname}\global\let\@Xout\@Xsub
            \openout\@Xout\@tempa.aux \W@X{\relax}%
            \DN@{\let\end\endinput\@input\@inputname
                    \@Xclose@\make@Xmain}\fi\fi\fi
  \donext@}
\def\includeonly#1{\edef\@includelist{#1}}
\def\noinclude#1{\edef\@noincllist{#1}}

%%%%%%%%%%%%%%%%%%%% Numeric styles
\def\arabicnum#1{\number#1}

\def\Romannum#1{\expandafter\uppercase\expandafter{\romannumeral#1}}
\def\alphnum#1{\ifcase#1\or a\or b\or c\or d\else\@ialph{#1}\fi}
\def\@ialph#1{\ifcase#1\or \or \or \or \or e\or f\or g\or h\or i\or j\or
    k\or l\or m\or n\or o\or p\or q\or r\or s\or t\or u\or v\or w\or x\or y\or
    z\else\fi}
\def\Alphnum#1{\ifcase#1\or A\or B\or C\or D\else\@Ialph{#1}\fi}
\def\@Ialph#1{\ifcase#1\or \or \or \or \or E\or F\or G\or H\or I\or J\or
    K\or L\or M\or N\or O\or P\or Q\or R\or S\or T\or U\or V\or W\or X\or Y\or
    Z\else\fi}

%%%%%%%%%%%%%%%%%%%% Counters
\def\ST@P{step}
\def\ST@LE{style}
\def\N@M{no}
\def\F@NT{font@}
%   #1 is the name of the counter to be defined,
%   #2 is the counter this one depends upon,
\outer\def\newcounter{\checkbrack@{\expandafter\newcounter@\@txtopt@{{}}}}
{\let\newcount\relax
\gdef\newcounter@#1#2#3{{%
    \toks@@\expandafter{\csname\eat@bs#2\N@M\endcsname}%
    \DN@{\alloc@0\count\countdef\insc@unt}%
    \ifx\@txtopt@\identity@\expandafter\next@\the\toks@@
        \else\if\notempty{#1}\global\@Nlet#2\N@M\@#1\fi\fi
    \@Nxdef\the\eat@bs#2\@{\if\@undefined\the\eat@bs#3\@\else
            \@Nx\the\eat@bs#3\@.\fi\noexpand\arabicnum\the\toks@@}%
  \@Nxdef#2\ST@P\@{}%
  \if\@undefined#3\ST@P\@\else
    \edef\next@{\noexpand\@@addto\@Nx#3\ST@P\@{%
             \global\@Nx#2\N@M\@\z@\@Nx#2\ST@P\@}}\next@\fi
    \expandafter\@@addto\expandafter\@Xclose\expandafter
        {\expandafter\\\the\toks@@}}}}
\outer\def\copycounter#1#2{%
    \@Nxdef#1\N@M\@{\@Nx#2\N@M\@}%
    \@Nxdef#1\ST@P\@{\@Nx#2\ST@P\@}%
    \@Nxdef\the\eat@bs#1\@{\@Nx\the\eat@bs#2\@}}
\outer\def\everystep{\checkstar@\everystep@}
\def\everystep@#1{\if@star\let\next@\gdef\else\let\next@\@@addto\fi
    \@DO\next@#1\ST@P\@}
%   #1 is the counter whose style is to be changed,
\def\counterstyle#1{\@Ngdef\the\eat@bs#1\@}
\def\advancecounter#1#2{\@N#1\ST@P\@\global\advance\@N#1\N@M\@#2}
\def\setcounter#1#2{\@N#1\ST@P\@\global\@N#1\N@M\@#2}
\def\counter#1{\refstepcounter#1\printcounter#1}
\def\printcounter#1{\@N\the\eat@bs#1\@}
\def\refcounter#1{\xdef\@lastmark{\printcounter#1}}
\def\stepcounter#1{\advancecounter#1\@ne}
\def\refstepcounter#1{\stepcounter#1\refcounter#1}
\def\savecounter#1{\@Nedef#1@sav\@{\global\@N#1\N@M\@\the\@N#1\N@M\@}}
\def\restorecounter#1{\@Name#1@sav\@}

%%%%%%%%%%%%%%%%%%%% Warnings
\def\warning#1#2{\W@{Warning: #1 on input line #2}}
\def\warning@#1{\warning{#1}{\the\inputlineno}}
\def\wrn@@Protect#1#2{\warning@{\string\Protect\string#1\space ignored}}
\def\wrn@@label#1#2{\warning{label `#1' multiply defined}{#2}}
\def\wrn@@ref#1#2{\warning@{label `#1' undefined}}
\def\wrn@@cite#1#2{\warning@{citation `#1' undefined}}
\def\wrn@@command#1#2{\warning@{Preamble command \string#1\space ignored}#2}
\def\wrn@@option#1#2{\warning@{Option \string#1\string#2\space is not supported}}
\def\wrn@@reference#1#2{\W@{Reference `#1' on input line \the\inputlineno}}
\def\wrn@@citation#1#2{\W@{Citation `#1' on input line \the\inputlineno}}
\let\wrn@reference\eat@ii
\let\wrn@citation\eat@ii
%% Disable wornings: works with
%% \Protect, \label, \ref, \cite, \command, \option, \reference, \citation
%% (with FONT.DEF)   \font
%% (with DEBUG.DEF)  \nocite
\def\nowarning#1{\if\@undefined\wrn@\eat@bs#1\@\wrn@option\nowarning#1\else
        \@Nlet\wrn@\eat@bs#1\@\eat@ii\fi}
\def\printwarning#1{\if\@undefined\wrn@@\eat@bs#1\@\wrn@option\printwarning#1\else
        \@Nlet\wrn@\eat@bs#1\expandafter\@\csname wrn@@\eat@bs#1\endcsname\fi}
\printwarning\Protect \printwarning\label \printwarning\ref
\printwarning\cite \printwarning\command \printwarning\option

%%%%%%%%%%%%%%%%%%%% Hyperrefs
{\catcode`\#=12\gdef\@lH{#}}
\def\@@HREF#1{}
\def\@HREF#1#2{\@@HREF{a #1}{\let\@@HREF\eat@#2}\@@HREF{/a}}
\def\@@Hf#1{file:#1} \let\@Hf\@@Hf
\def\@@Hl#1{\if\notempty{#1}\@lH#1\fi} \let\@Hl\@@Hl
\def\@@Hname#1{\@HREF{name="#1"}{}} \let\@Hname\@@Hname
\def\@@Href#1{\@HREF{href="#1"}} \let\@Href\@@Href
\ifx\undefined\pdfoutput
  \csname newcount\endcsname\pdfoutput
\else
  \def\pdflinkattr{attr{/C [0 0.9 0.9]}}
  \let\pdflinkbegin\empty
  \let\pdflinkend\empty
  \def\@pdfHf#1{file {#1}}
  \def\@pdfHl#1{name {#1}}
  \def\@pdfHname#1{\pdfdest name{#1}xyz\relax}
  \def\@pdfHref#1#2{\pdfstartlink \pdflinkattr goto #1\relax
    \pdflinkbegin#2\pdflinkend\pdfendlink}
  \def\@ifpdf#1#2{\ifnum\pdfoutput>\z@\expandafter#1\else\expandafter#2\fi}
  \def\@Hf{\@ifpdf\@pdfHf\@@Hf}
  \def\@Hl{\@ifpdf\@pdfHl\@@Hl}
  \def\@Hname{\@ifpdf\@pdfHname\@@Hname}
  \def\@Href{\@ifpdf\@pdfHref\@@Href}
\fi
\def\@Hr#1#2{\if\notempty{#1}\@Hf{#1}\fi\@Hl{#2}}
\def\@localHref#1{\@Href{\@Hr{}{#1}}}
\def\@countlast#1{\@N#1last\@}
\def\@@countref#1#2{\global\advance#2\@ne
  \@Nxdef#2last\@{\the#2}\@tocHname{#1\@countlast#2}}
\def\@countref#1{\@DO\@@countref#1@HR\@#1}

%\Href@-xxx#1#2#3
\def\Href@@#1{\@N\Href@-#1\@}
\def\Href@#1#2{\@N\Href@-#1\@{\@Hl{@#1-#2}}}
%\Hname@-xxx#1
\def\Hname@#1{\@N\Hname@-#1\@}
%\Hlast@-xxx
\def\Hlast@#1{\@N\Hlast@-#1\@}
\def\cntref@#1{\global\@DO\advance\cnt@#1\@\@ne
  \@Nxdef\Hlast@-#1\@{\@DO\the\cnt@#1\@}\Hname@{#1}{@#1-\Hlast@{#1}}}
\def\HyperRefs#1{\global\@Nlet\Hlast@-#1\@\empty
  \global\@Nlet\Hname@-#1\@\@Hname
  \global\@Nlet\Href@-#1\@\@Href}
\def\NoHyperRefs#1{\global\@Nlet\Hlast@-#1\@\empty
  \global\@Nlet\Hname@-#1\@\eat@
  \global\@Nlet\Href@-#1\@\eat@}

%% Hyperrefs for label/ref
\HyperRefs{label} {\catcode`\-11
\gdef\@labelref#1{\Hname@-label{r@-#1}}
\gdef\@xHref#1{\Href@-label{\@Hl{r@-#1}}} }
%% Hyperrefs for toc
\HyperRefs{toc}
\def\@HR#1{\if\notempty{#1}\string\@HR{\Hlast@{toc}}{#1}\else{}\fi}

%\def\HyperRefs{\gdef\@labelref##1{\@Hname{r@-##1}}%
% \gdef\@xHref##1{\@localHref{r@-##1}}}
%\def\NoHyperRefs{\global\let\@labelref\eat@\global\let\@xHref\eat@}
%\def\TOCRefs{\global\let\@tocHname\@Hname\global\let\@tocHref\@localHref}
%\def\NoTOCRefs{\global\let\@tocHname\eat@\global\let\@tocHref\eat@}
%%\def\@HR#1{\if\notempty{#1}\string\@HR{\cnt@toclast}{#1}\else{}\fi}
%\HyperRefs
%%\TOCRefs

%%\def\@CLR#1{\special{color push #1}}
%%\def\@eCLR{\special{color pop}}

%%%%%%%%%%%%%%%%%%%% Labels, x-references
\def\bftext{\ifmmode\fam\bffam\else\bf\fi}
\let\@lastmark\empty
\let\@lastlabel\empty
\def\lastmark{\@lastmark}
\let\lastlabel\empty
\let\everylabel\relax
\let\everylabel@\eat@
\let\everyref\relax
\def\newlabel{\bgroup\everylabel\newlabel@}
\def\newlabel@#1#2#3{\if\@undefined\r@-#1\@\else\wrn@label{#1}{#3}\fi
  {\let\protect\noexpand\@Nxdef\r@-#1\@{#2}}\egroup}
\def\w@ref{\bgroup\everyref\w@@ref}
\def\w@@ref#1#2#3#4{%
  \if\@undefined\r@-#1\@{\bftext??}#2{#1}{}\else%
   \@xHref{#1}{\@DO{\expandafter\expandafter#3}\r@-#1\@\@nil}\fi
  #4{#1}{}\egroup}%\null}}
\def\@@@xref#1{\w@ref{#1}\wrn@ref\@car\wrn@reference}
\def\@xref#1{\rom{\@@@xref{#1}}}
\let\xref\@xref
\def\pageref#1{\w@ref{#1}\wrn@ref\@cdr\wrn@reference}
\def\thepage{\ifnum\pageno<\z@\romannumeral-\pageno\else\number\pageno\fi}
\def\label@{\@bsphack\bgroup\everylabel\label@@}
\def\label@@#1#2{\everylabel@{{#1}{#2}}%
  \@labelref{#2}%
  \let\thepage\relax
  \def\protect{\noexpand\noexpand\noexpand}%
  \edef\@tempa{\edef\noexpand\@lastlabel{#1}%
    \W@X{\string\newlabel{#2}{{\@lastmark}{\thepage}}{\the\inputlineno}}}%
  \expandafter\egroup\@tempa\@esphack}
\def\label#1{\label@{#1}{#1}}
\def\fn@P@{\relaxnext@
    \DN@{\ifx[\next\DN@[####1]{}\else
        \ifx"\next\DN@"####1"{}\else\DN@{}\fi\fi\next@}%
    \FN@\next@}
\def\eat@fn#1{\ifx#1[\expandafter\eat@br\else
  \ifx#1"\expandafter\expandafter\expandafter\eat@qu\fi\fi}
\def\eat@br#1]#2{}
\def\eat@qu#1"#2{}
{\catcode`\~\active\lccode`\~`\@
\lowercase{\global\let\@@P@~\gdef~{\protect\@@P@}}}
\def\Protect@@#1{\def#1{\protect#1}}
\def\disable@special{\let\W@X@\eat@iii\let\label\eat@
    \def\footnotemark{\protect\fn@P@}%
  \let\footnotetext\eat@fn\let\footnote\eat@fn
    \let\refcounter\eat@\let\savecounter\eat@\let\restorecounter\eat@
    \let\advancecounter\eat@ii\let\setcounter\eat@ii
  \let\ifvmode\iffalse\Protect@@\@@@xref\Protect@@\pageref\Protect@@\nofrills
    \Protect@@\\\Protect@@~}
\let\notoctext\identity@
\def\W@X@#1#2#3{\@bsphack{\disable@special\let\notoctext\eat@
    \def\chapter{\protect\chapter@toc}\let\thepage\relax
    \def\protect{\noexpand\noexpand\noexpand}#1%
  \edef\next@{\if\@undefined#2\@\else\write#2{#3}\fi}\expandafter}\next@
    \@esphack}
\newcount\cnt@toc
\def\writeauxline#1#2#3{\W@X@{\cntref@{toc}\let\tocref\@HR}
  \@Xout{\string\@Xline{#1}{#2}{#3}{\thepage}}}
{\let\newwrite\relax
\gdef\@openin#1{\make@letter\@input{\jobname.#1}\t@mpcat}
\gdef\@openout#1{\global\expandafter\newwrite\csname
tf@-#1\endcsname
   \immediate\openout\@N\tf@-#1\@\jobname.#1\relax}}
\def\@@openout#1{\@openout{#1}%
  \@@addto\readaux@{\immediate\closeout\@N\tf@-#1\@}}
\def\auxlinedef#1{\@Ndef\do@-#1\@}
\def\@Xline#1{\if\@undefined\do@-#1\@\expandafter\eat@iii\else
    \@DO\expandafter\do@-#1\@\fi}
\def\beginW@{\bgroup\def\do##1{\catcode`##112 }\dospecials\do\@\do\"
    \catcode`\{\@ne\catcode`\}\tw@\immediate\write\@N}
\def\endW@toc#1#2#3{{\string\tocline{#1}{#2\string\page{#3}}}\egroup}
\def\do@tocline#1{%
%%  The file version
    \if\@undefined\tf@-#1\@\expandafter\eat@iii\else
        \beginW@\tf@-#1\@\expandafter\endW@toc\fi
%%  The \toks version
%       \if\@undefined\the#1@@\@\else
%           \global\addto{\the#1@}{\tocline{#2}{#3\page{#4}}}\fi
} \auxlinedef{toc}{\do@tocline{toc}}

\let\protect\empty
\def\Protect#1{\if\@undefined#1@P@\@\PROTECT#1\else\wrn@Protect#1\empty\fi}
\def\PROTECT#1{\@Nlet#1@P@\@#1\edef#1{\noexpand\protect\@Nx#1@P@\@}}
\def\pdef#1{\edef#1{\noexpand\protect\@Nx#1@P@\@}\@Ndef#1@P@\@}

\Protect\operatorname \Protect\operatornamewithlimits
\Protect\qopname@ \Protect\qopnamewl@ \Protect\text
\Protect\topsmash \Protect\botsmash \Protect\smash
\Protect\widetilde \Protect\widehat \Protect\thetag
\Protect\therosteritem
% Fonts:
\Protect\Cal \Protect\Bbb \Protect\bold \Protect\slanted
\Protect\roman \Protect\italic \Protect\boldkey
\Protect\boldsymbol \Protect\frak \Protect\goth \Protect\dots
% Symbols
\Protect\cong \Protect\lbrace \let\{\lbrace \Protect\rbrace
\let\}\rbrace
\let\root@P@@\root \def\root@P@#1{\root@P@@#1\of}
\def\root#1\of{\protect\root@P@{#1}}

\def\frills{\ignorespaces\@txtopt@}
\def\frillsnotempty#1{\x@notempty{\@txtopt@{#1}}}
\def\numberline{\@numopt@}
\newif\if@theorem
\let\@therstyle\eat@
\def\@headtext@#1#2{{\disable@special\let\protect\noexpand
    \def\chapter{\protect\chapter@rh}%
    \edef\next@{\noexpand\@frills@\noexpand#1{#2}}\expandafter}\next@}
\let\AmSrighthead@\rightheadtext
\def\rightheadtext{\checkfrills@{\@headtext@\AmSrighthead@}}
\let\AmSlefthead@\leftheadtext
\def\leftheadtext{\checkfrills@{\@headtext@\AmSlefthead@}}
% #1 refers to the style,
% #2 refers to the style in the toc,
% #3 refers to the counter,
% #4 represents the AmSTeX's counterpart (we use \end... because
%   of \outer-ness), and
% #5 is the text to be typeset
\def\@head@@#1#2#3#4#5{\@Name\pre\eat@bs#1\@\if@theorem\else
    \@frills@{\csname\expandafter\eat@iv\string#4\endcsname}\relax
        \ifx\protect\empty\@N#1\F@NT\@\fi\fi
    \@N#1\ST@LE\@{\counter#3}{#5}%
  \if@write\writeauxline{toc}{\eat@bs#1}{#2{\counter#3}\@HR{#5}}\fi
    \if@theorem\else\expandafter#4\fi
    \ifx#4\endhead\ifx\@txtopt@\identity@\else
        \headmark{\@N#1\ST@LE\@{\counter#3}{\frills\empty}}\fi\fi
    \@Name\post\eat@bs#1\@\ignorespaces}
\ifx\undefined\endhead\Invalid@\endhead\fi
\def\@head@#1{\checkstar@{\checkfrills@{\checkbrack@{\@head@@#1}}}}
% #1 is the name,
% #2 is the counter, and
% #3 is the title text
\def\@thm@@#1#2#3{\@Name\pre\eat@bs#1\@
    \@frills@{\csname\expandafter\eat@iv\string#3\endcsname}
    {\@theoremtrue\check@therstyle{\@N#1\ST@LE\@}\frills
            {\counter#2}\@theoremfalse}%
    \@DO\envir@stack\end\eat@bs#1\@
    \@N#1\F@NT\@\@Name\post\eat@bs#1\@\ignorespaces}
\def\@thm@#1{\checkstar@{\checkfrills@{\checkbrack@{\@thm@@#1}}}}
% #1 is the name,
% #2 is the counter,
% #3 is the name of the corresponding table (lof, lot, etc.)
% #4 is either \topcaption or \botcaption
%   #5 is the caption text
\def\@capt@@#1#2#3#4#5\endcaption{\bgroup
    \edef\@tempb{\global\footmarkcount@\the\footmarkcount@
    \global\@N#2\N@M\@\the\@N#2\N@M\@}%
    \def\shortcaption##1{\global\def\sh@rtt@xt####1{##1}}\let\sh@rtt@xt\identity@
%    \let\notoctext\identity@
%%% To handle \nofrills !!!
%   \DN@##1##2##3{\false@\fi\iftrue}%
%   \ifx\@frills@\identity@\else\let\notempty\next@\fi
    \DN@{#4{\@tempb\@N#1\ST@LE\@{\counter#2}}}%
    \if\notempty{#5}\DNii@{\next@\@N#1\F@NT\@}\else\let\nextii@\next@\fi
    \nextii@#5\endcaption
  \if@write\writeauxline{#3}{\eat@bs#1}{{} \@HR{\@N#1\ST@LE\@{\counter#2}%
    \if\notempty{#5}.\enspace\fi\sh@rtt@xt{#5}}}\fi
  \global\let\sh@rtt@xt\undefined\egroup}
\def\@capt@#1{\checkstar@{\checkfrills@{\checkbrack@{\@capt@@#1}}}}
\let\captiontextfont@\empty

\ifx\undefined\subsubheadfont@\def\subsubheadfont@{\it}\fi
\ifx\undefined\proclaimfont\def\proclaimfont{\sl}\fi
\ifx\undefined\proclaimfont@\let\proclaimfont@\proclaimfont\fi
\def\proclaimfont{\proclaimfont@}
\ifx\undefined\definitionfont@\def\AmSdeffont@{\rm}
    \else\let\AmSdeffont@\definitionfont@\fi
\ifx\undefined\remarkfont@\def\remarkfont@{\rm}\fi

\def\newfont@def#1#2{\if\@undefined#1\F@NT\@
    \@Nxdef#1\F@NT\@{\@Nx.\expandafter\eat@iv\string#2\F@NT\@}\fi}
% #1 is the name (and the style),
% #2 is the style in the toc,
% #3 is the counter, and
% #4 is the AmSTeX's counterpart to be used (in the \end... form):
\def\newhead@#1#2#3#4{{%
    \gdef#1{\@therstyle\@therstyle\@head@{#1#2#3#4}}\newfont@def#1#4%
    \if\@undefined#1\ST@LE\@\@Ngdef#1\ST@LE\@{\headstyle}\fi
    \if\@undefined#2\@\gdef#2{\headtocstyle}\fi
  \@@addto\moretocdefs@{\\#1#1#4}}}
\outer\def\newhead#1{\checkbrack@{\expandafter\newhead@\expandafter
    #1\@txtopt@\headtocstyle}}
% #1 is the default title (like Theorem, Lemma, etc.),
% #2 is the name to be defined,
% #3 refers to the counter (which should be defined separately),
% #4 is the AmSTeX's counterpart: \endproclaim, \endremark, \endAmSdef
\outer\def\newtheorem#1#2#3#4{{%
    \gdef#2{\@thm@{#2#3#4}}\newfont@def#2#4%
    \@Nxdef\end\eat@bs#2\@{\noexpand\revert@envir
        \@Nx\end\eat@bs#2\@\noexpand#4}%
  \if\@undefined#2\ST@LE\@\@Ngdef#2\ST@LE\@{\proclaimstyle{#1}}\fi}}%
% #1 is the default title (like Figure, Table, etc.),
% #2 is the name to be defined,
% #3 refers to the counter (which should be defined separately),
% #4 is the name of the corresponding table (toc, lof, lot, etc.)
% #5 is either \topcaption or \botcaption
\outer\def\newcaption#1#2#3#4#5{{\let#2\relax
  \edef\@tempa{\gdef#2####1\@Nx\end\eat@bs#2\@}%
    \@tempa{\@capt@{#2#3{#4}#5}##1\endcaption}\newfont@def#2\endcaptiontext%
  \if\@undefined#2\ST@LE\@\@Ngdef#2\ST@LE\@{\captionstyle{#1}}\fi
  \@@addto\moretocdefs@{\\#2#2\endcaption}\newtoc{#4}}}
{
% #1 is the name of the table (toc, lof, lot, etc.)
\outer\gdef\newtoc#1{{%
    \@DO\ifx\do@-#1\@\relax
%%  The \toks version
%    \global\expandafter\newtoks\csname the#4@@\endcsname
    \global\auxlinedef{#1}{\do@tocline{#1}}{}%
    \@@addto\tocsections@{\make@toc{#1}{}}\fi}}}

\toks@\expandafter{\itembox@}
\toks@@{\bgroup\let\therosteritem\identity@\let\rm\empty
  \let\@Href\eat@\let\@Hname\eat@
  \edef\next@{\edef\noexpand\@lastmark{\therosteritem@}}\donext@}
\edef\itembox@{\the\toks@@\the\toks@}
\def\firstitem@false{\let\iffirstitem@\iffalse
    \global\let\lastlabel\@lastlabel}

\let\rosteritemrefform\therosteritem
\let\rosteritemrefseparator\empty
\def\rosteritemref#1{\hbox{\rosteritemrefform{\@@@xref{#1}}}}
\def\local#1{\label@\@lastlabel{\lastlabel-i#1}}

\def\xRef@P@{\gdef\lastlabel}
\def\xRef#1{\@xref{#1}\protect\xRef@P@{#1}}

\def\iref@P@{\gdef\lastref}
\def\itemref#1#2{\rosteritemref{#1-i#2}\protect\iref@P@{#1}}
\def\iref#1{\@xref{#1}\rosteritemrefseparator\itemref{#1}}

\def\eqref#1{\thetag{\@@@xref{#1}}}
\def\tagform@#1{\ifmmode\hbox{\rm\else\rom{\fi
        (\ignorespaces#1\unskip)\iftrue}\else}\fi}

\let\AmSfnote@\makefootnote@
\def\makefootnote@#1{\bgroup\let\footmarkform@\identity@
  \edef\next@{\edef\noexpand\@lastmark{#1}}\donext@\AmSfnote@{#1}}

\def\clearpage{\ifnum\insertpenalties>0\line{}\fi\vfill\supereject}

\def\proof{\checkfrills@{\checkbrack@{%
    \check@therstyle{\@frills@{\demo}{\frills{Proof}}{}}
        {\frills{}\envir@stack\endremark\envir@stack\enddemo}%
  \envir@stack\endproof\ignorespaces}}}
\def\everyendproof{\qed}
\def\endproof{\nofrillscheck{\frills@{\everyendproof}\revert@envir\endproof\enddemo}}

\let\AmSref\ref
\let\AmSrefstyle\refstyle
\let\plaincite\cite
\def\citei@#1,{\citeii@#1\eat@,}
\def\citeii@#1\eat@{\w@ref{#1}\wrn@cite\@car\wrn@citation}
\def\mcite@#1;{\plaincite{\citei@#1\eat@,\unskip}\mcite@i}
\def\mcite@i#1;{\DN@{#1}\ifx\next@\endmcite@
  \else, \plaincite{\citei@#1\eat@,\unskip}\expandafter\mcite@i\fi}
\def\endmcite@{\endmcite@}
\def\cite#1{\mcite@#1;\endmcite@;}
\PROTECT\cite
\def\refstyle#1{\AmSrefstyle{#1}\uppercase{%
    \ifx#1A\relax \def\@ref@##1{\AmSref\xdef\@lastmark{##1}\key##1}%
    \else\ifx#1C\relax \def\@ref@##1{\AmSref\no\counter\refno}%
        \else\def\@ref@{\AmSref}\fi\fi}}
\refstyle A
\newcounter\refno\null
\newif\ifRefs
\gdef\Refs{\checkstar@{\checkbrack@{\csname AmSRefs\endcsname
  \nofrills{\frills{References}%
  \if@write\writeauxline{toc}{vartocline}{\@HR{\frills{References}}}\fi}%
  \def\ref{\@ref@}\Refstrue\ignorespaces}}}
\let\ref\xref

\newif\iftoc
\pdef\tocbreak{\iftoc\hfil\break\fi}
\def\tocsections@{\make@toc{toc}{}}
\let\moretocdefs@\empty
\def\newtocline@#1#2#3{%
  \edef#1{\def\@Nx#2line\@####1{\@Nx.\expandafter\eat@iv
        \string#3\@####1\noexpand#3}}%
  \@Nedef\no\eat@bs#1\@{\let\@Nx#2line\@\noexpand\eat@}%
    \@N\no\eat@bs#1\@}
\def\MakeToc#1{\@@openout{#1}}
\def\newtocline#1#2#3{\Err@{\Invalid@@\string\newtocline}}
\def\make@toc#1#2{\penaltyandskip@{-200}\aboveheadskip
    \if\notempty{#2}
        \centerline{\headfont@\ignorespaces#2\unskip}\nobreak
    \vskip\belowheadskip \fi
%%  The file version
    \@openin{#1}\relax
%%  The \toks version
%   \if\@undefined\the#1@@\@
%       \the\@N\the#1@@\@\global\@N\the#1@@\@\fi
    \vskip\z@}
\def\contents{\readaux\checkfrills@{\checkbrack@{\@contents@}}}
\def\@contents@{\toc@{\frills{Contents}}\envir@stack\endcontents%
    \def\nopagenumbers{\let\page\eat@}\let\newtocline\newtocline@\toctrue
  \def\@HR{\Href@{toc}}%
  \def\tocline##1{\csname##1line\endcsname}
  \edef\caption##1\endcaption{\expandafter\noexpand
    \csname head\endcsname##1\noexpand\endhead}%
    \ifmonograph@\def\vartoclineline{\Chapterline}%
        \else\def\vartoclineline##1{\sectionline{{} ##1}}\fi
  \let\\\newtocline@\moretocdefs@
    \ifx\@frills@\identity@\def\\##1##2##3{##1}\moretocdefs@
        \else\let\tocsections@\relax\fi
    \def\\{\unskip\space\ignorespaces}\let\maketoc\make@toc}
\def\endcontents{\tocsections@\vskip-\lastskip\revert@envir\endcontents
    \endtoc}

% \selectf@nt is for future extensions (like RUSSIAN.TEX)
\if\@undefined\selectf@nt\@\let\selectf@nt\identity@\fi
\def\Err@math#1{\Err@{Use \string#1\space only in text}}
\def\textonlyfont@#1#2{%
    \def#1{\RIfM@\Err@math#1\else\edef\f@ntsh@pe{\string#1}\selectf@nt#2\fi}%
    \PROTECT#1}
\tenpoint

% #1 is the name of the switch to be defined
%   #2 is the default font switch
\def\newshapeswitch#1#2{\gdef#1{\selectsh@pe#1#2}\PROTECT#1}
% #1 is the name of the switch
% #2 is the current shape
% #3 is the shape to be used
\def\shapeswitch#1#2#3{\@Ngdef#1\string#2\@{#3}}
% These shapes are used by \rom
\shapeswitch\rm\bf\bf \shapeswitch\rm\tt\tt
\shapeswitch\rm\smc\smc
\newshapeswitch\em\it
% These shapes are used by \em and \emph
\shapeswitch\em\it\rm \shapeswitch\em\sl\rm
\def\selectsh@pe#1#2{\relax\if\@undefined#1\f@ntsh@pe\@#2\else
    \@N#1\f@ntsh@pe\@\fi}

\def\@itcorr@{\leavevmode
    \edef\prevskip@{\ifdim\lastskip=\z@ \else\hskip\the\lastskip\relax\fi}\unskip
    \edef\prevpenalty@{\ifnum\lastpenalty=\z@ \else
        \penalty\the\lastpenalty\relax\fi}\unpenalty
    \/\prevpenalty@\prevskip@}
\def\rom@P@#1{\@itcorr@{\selectsh@pe\rm\rm#1}}
\def\rom{\protect\rom@P@}
%%%%    \Rom will unconditionally switch to \rm, like \rom in AmSppt
\def\Rom@P@#1{\@itcorr@{\rm#1}}
\def\Rom{\protect\Rom@P@}
{\catcode`\-11 \HyperRefs{idx} \HyperRefs{glo}
\newcount\cnt@idx \global\cnt@idx=10000
\newcount\cnt@glo \global\cnt@glo=10000
\gdef\writeindex#1{\W@X@{\cntref@{idx}}\tf@-idx
 {\string\indexentry{#1}{\Hlast@{idx}}{\thepage}}}
\gdef\writeglossary#1{\W@X@{\cntref@{glo}}\tf@-glo
 {\string\glossaryentry{#1}{\Hlast@{glo}}{\thepage}}}
}
\def\emph#1{\@itcorr@\bgroup\em\ignorespaces#1\unskip\egroup
  \DN@{\DN@{}\ifx\next.\else\ifx\next,\else\DN@{\/}\fi\fi\next@}\FN@\next@}
\def\makequoteactive{\catcode`\"\active}
{\makequoteactive\gdef"{\FN@\quote@}
\gdef\quote@{\ifx"\next\DN@"##1""{\quoteii{##1}}\else\DN@##1"{\quotei{##1}}\fi\next@}}
\let\quotei\eat@
\let\quoteii\eat@
\def\MakeIndex{\@openout{idx}}
\def\MakeGlossary{\@openout{glo}}

%%%%%%%%%%%%%%%%%%%%%%  Just helpful things %%%%%%%%%%%%%%%%%%%%%%%%%%%%
\def\endofpar#1{\ifmmode\ifinner\endofpar@{#1}\else\eqno{#1}\fi
    \else\leavevmode\endofpar@{#1}\fi}
\def\endofpar@#1{\unskip\penalty\z@\null\hfil\hbox{#1}\hfilneg\penalty\@M}

\newdimen\normalparindent\normalparindent\parindent
\def\firstparindent#1{\everypar\expandafter{\the\everypar
  \global\parindent\normalparindent\global\everypar{}}\parindent#1\relax}

%% Commands to disable
\@@addto\disablepreambule@cs{%
    \\\readaux\relax
    \\\begin\relax
    \\\readaux@\relax
    \\\@openout\eat@
    \\\@@openout\eat@
    \/\Monograph\empty
    \/\MakeIndex\empty
    \/\MakeGlossary\empty
    \/\MakeToc\eat@
    \/\HyperRefs\eat@
    \/\NoHyperRefs\eat@
}

\csname label.def\endcsname

%%%%%%%%%%%%%%%%%%%%%%%% Definitions for my papers %%%%%%%%%%%%%%%%%%%%%%%

\def\punct#1#2{\if\notempty{#2}#1\fi}
\def\sppunct{\punct{.\enspace}}
\def\varpunct#1#2{\if\frillsnotempty{#2}#1\fi}

\def\headstyle#1#2{\numberline{#1\sppunct{#2}}\ignorespaces#2\unskip}
\def\headtocstyle#1#2{\numberline{#1\punct.{#2}}\space #2}

\def\specialtocstyle#1#2{#2}
\newcounter\section\null
\newcounter\subsection\section
\newcounter\subsubsection\subsection
\newhead\specialsection[\specialtocstyle]\null\endspecialhead
\newhead\section\section\endhead
\newhead\subsection\subsection\endsubhead
\newhead\subsubsection\subsubsection\endsubsubhead
\def\firstappendix{\global\sectionno0 %
  \counterstyle\section{\Alphnum\sectionno}%
    \global\let\firstappendix\empty}

\def\appendixtocstyle#1#2{\space\numberline{Appendix #1\sppunct{#2}}#2}
\newhead\appendix[\appendixtocstyle]\section\endhead

\let\endAmSdef\enddefinition
\def\proclaimstyle#1#2{\numberline{#2\varpunct{.\enspace}{#1}}\frills{#1}}
\copycounter\thm\subsubsection
%\newcounter[\subsubsectionno]\thm\subsection
\theorem\thm\endproclaim
\proposition\thm\endproclaim
\lemma\thm\endproclaim
\corollary\thm\endproclaim
\definition\thm\endAmSdef
\example\thm\endAmSdef

\def\captionstyle#1#2{\frills{#1}\numberline{\varpunct{ }{#1}#2}}
\newcounter\figure\null
\newcounter\table\null
\newcaption{Figure}\figure\figure{lof}\botcaption
\newcaption{Table}\table\table{lot}\topcaption

\copycounter\equation\subsubsection
%\newcounter[\subsubsectionno]\equation\subsection

%\endinput

\expandafter\ifx\csname label.def\endcsname\relax\input label.def
\fi
\def\stydate{June 26, 2000}
\def\styname{DEBUG.DEF}
\immediate\write16{This is \styname\space by A.Degtyarev
<\stydate>} \expandafter\ifx\csname debug.def\endcsname\relax\else
  \message{[already loaded]}\endinput\fi
\expandafter\edef\csname debug.def\endcsname{%
  \catcode`\noexpand\@\the\catcode`\@\edef\noexpand\styname{\styname}
  \def\expandafter\noexpand\csname debug.def\endcsname{\stydate}}
\catcode`\@=11 {\edef\temp{\the\everyjob\W@{\styname: <\stydate>}}
\global\everyjob\expandafter{\temp}}

%%%%%%%%%%%%%%%%%%%% Marginal notes
\def\n@te#1#2{\leavevmode\vadjust{%
 {\setbox\z@\hbox to\z@{\strut\eightpoint\let\quotei\filename#1}%
  \setbox\z@\hbox{\raise\dp\strutbox\box\z@}\ht\z@=\z@\dp\z@=\z@%
  #2\box\z@}}}
\def\leftnote#1{\n@te{\hss#1\quad}{}}
\def\rightnote#1{\n@te{\quad\kern-\leftskip#1\hss}{\moveright\hsize}}
\def\?{\FN@\qumark}
\def\qumark{\ifx\next"\DN@"##1"{\leftnote{\rm##1}}\else
 \DN@{\leftnote{\rm??}}\fi{\rm??}\next@}
\def\filename#1{\hbox{\tt #1}}
\def\mnote@@#1{\rightnote{\vtop{%
 \ifcat\noexpand"\noexpand~\def"##1"{\filename{##1}}\fi
 \hsize2.0in \baselineskip7\p@\parindent\z@
 \tolerance\@M\spaceskip2.6\p@ plus10\p@ minus.9\p@\rm#1}}}
\def\mnote#1{\@bsphack\mnote@{#1}\@esphack}

\def\nonotes{\let\mnote@\eat@}
\def\printnotes{\let\mnote@\mnote@@}
\printnotes

%%%%%%%%%%%%%%%%%%%% Printing labels
\def\PrintLabels{%
 \gdef\printlabel@##1##2{\ifvmode\else\leftnote{\eighttt##2}\fi}}
\def\NoLabels{\global\let\printlabel@\eat@ii}
\NoLabels \@@addparm\everylabel@{\printlabel@#1}

%%%%%%%%%%%%%%%%%%%% Printing file info
\def\PrintFiles{\gdef\outputmark@{\line{\hfill\smash{\raise1cm\vbox{%
  \hbox to\z@{\kern1cm\tenrm\the\month/\the\day/\the\year\hss}%
  \hbox to\z@{\kern1cm\tentt\jobname.tex\hss}%
  \ifx\filecomment\undefined\else\hbox to\z@{\kern1cm\tenrm\filecomment\hss}\fi}}}}}
\def\NoFiles{\global\let\outputmark@\empty}
\def\PageMark{\gdef\outputmark@}
\let\outputmark@\empty
\@@addto\pagetop@\outputmark@

%%%%%%%%%%%%%%%%%%%% Tracing references
%\def\tracerefs{\def\wrn@reference##1##2{\@Nxdef\ref@-##1\@{##1}%
%    \W@X{\string\reference{##1}}}}
\def\tracerefs{\def\wrn@reference##1##2{\W@X{\string\reference{##1}}}}
\def\tracecites{\printwarning\nocite
  \def\wrn@citation##1##2{\@Nxdef\cite@-##1\@{##1}\W@X{\string\citation{##1}}}}
\def\wrn@@nocite#1#2{\ifRefs\wrn@@@nocite{#2}\fi}
\def\wrn@@@nocite#1{\if\@undefined\cite@-#1\@\warning@{citation `#1' [\@lastmark] not used}\fi}
\let\wrn@nocite\eat@ii
\let\reference\eat@
\let\citation\eat@
\@@addparm\everylabel@{\wrn@nocite#1}

\@@addto\disablepreambule@cs{%
    \/\PrintFiles\empty
    \/\NoFiles\empty
    \/\PageMark\empty
}

\csname debug.def\endcsname
%\endinput

\def\Per{\mathop{\roman{Per}}}
\def\conj{\operatorname{conj}}
\def\C{{\Bbb C}}
\def\R{{\Bbb R}}
\def\Z{{\Bbb Z}}
\def\Q{{\Bbb Q}}
\def\Rp#1{\R\roman P^{#1}}

\def\Im{\mathop{\roman{Im}}\nolimits}
\def\Re{\mathop{\roman{Re}}\nolimits}

\def\la{\langle}
\def\ra{\rangle}
\def\id{\mathop{\roman{id}}\nolimits}
\def\Aut{\mathop{\roman{Aut}}\nolimits}
\def\emptyset{\varnothing}
\def\oo{\varnothing}
\def\discr{\operatorname{discr}}
\def\Pic{\operatorname{Pic}}
\def\M{\Bbb M}
\def\MC{\M_0\otimes\C}
\def\PMC{P(\M_0\otimes\C)}
\def\D{\Delta}
\def\G{\Gamma}

\def\per{\mathop{\roman{per}}}
\def\Res{\mathop{\roman{Res}}}
\def\d{\roman{d}}

\def\til{\widetilde}

\def\L{\Lambda}
\def\s{\sigma}
\def\LL{\Bbb L}

\let\tm\proclaim
\let\endtm\endproclaim
\let\rk=\remark
\let\endrk=\endremark
\let\ge\geqslant
\let\le\leqslant
\let\+\sqcup
\let\dsum\+

\NoBlackBoxes \rightheadtext{On the deformation chirality of real
cubic fourfolds}

 \topmatter
\title
On the deformation chirality of real cubic fourfolds
\endtitle
\author S.~Finashin, V.~Kharlamov
\endauthor
\address Middle East Technical University,
Department of Mathematics\endgraf Ankara 06531 Turkey
\endaddress
%\email  serge@metu.edu.tr \endemail
\address
Universit\'{e} Louis Pasteur et IRMA (CNRS)\endgraf 7 rue Ren\'{e}
Descartes 67084 Strasbourg Cedex, France
\endaddress
%\affil \endaffil
%\address \endaddress
%\email \endemail
%\date \enddate
%\thanks \endthanks
%\dedicatory \enddedicatory
%\keywords \endkeywords
\abstract According to our previous results, the conjugacy class
of the involution induced by the complex conjugation in the
homology of a real non-singular cubic fourfold determines the
fourfold up to projective equivalence and deformation. Here, we
show how to eliminate the projective equivalence and to obtain a
pure deformation classification, that is how to respond to the
chirality question: which cubics are not deformation equivalent to
their image under a mirror reflection. We provide an arithmetical
criterion of chirality, in terms of the eigen-sublattices of the
complex conjugation involution in homology, and show how this
criterion can be effectively applied taking as examples $M$-cubics
(that is those for which the real locus has the richest topology)
and $(M-1)$-cubics (the next case with respect to complexity of
the real locus). It happens that there is one chiral class of
$M$-cubics and three chiral classes of $(M-1)$-cubics, contrary to
two achiral classes of $M$-cubics and three achiral classes of
$(M-1)$-cubics.
\endabstract
\endtopmatter

\document
\rightline{\vbox{\hsize 60mm \noindent\eightit\baselineskip10pt
L'univers est un ensemble dissym\'{e}trique, et je suis persuad\'{e}
que la vie, telle qu'elle manifeste \`{a} nous, est fonction de la
dissym\'{e}trie de l'univers ou des cons\'{e}quences qu'elle
entra\^{\i}ne. \vskip3mm Louis Pasteur
 \vskip3mm\noindent\eightrm
Observations sur les forces dissym\'{e}triques, CRAS, 78 (1874),
1515--1518
 \vskip10mm
}}

\heading \S1. Introduction \endheading

Recall that the projective nonsingular cubic fourfolds form the
complement in a projective space $P_{4,3}=P(Sym^3(\C^6))$ of
dimension $\binom{5+3}3-1=55$ to the so-called \emph{discriminant
hypersurface.} The discriminant hypersurface, which we denote by
$\Delta_{4,3}$, is defined over reals and its real part
$\Delta_{4,3}(\R)$ is represented by real singular cubics, so that
the space under our study is nothing but
$P_{4,3}(\R)\setminus\Delta_{4,3}(\R)$. (Such a notation specifies
the dimension, $n=4$, and the degree, $d=3$, of the hypersurfaces
under consideration; we make use of it in Section 8 for arbitrary
$n$ and $d$).

The space $\Cal{C}=P_{4,3}(\C)\setminus \Delta_{4,3}(\C)$ is
connected, while $\Cal{C}_\R=P_{4,3}(\R)\setminus\Delta_{4,3}(\R)$
is not. Understanding the nature of the connected components of
the latter space is a natural, and classical task, it can be
rephrased as a \emph{deformation classification} of real
projective nonsingular cubic fourfolds. In our previous paper
\cite{FK} we performed a classification with respect to a weaker
\emph{coarse deformation equivalence}: we call two real projective
nonsingular hypersurfaces {\it coarse deformation equivalent} if
one hypersurface is deformation equivalent to a projective
transformation of the other.

The difference between these two equivalence relations shows up in
the case of subvarieties of real projective spaces of odd
dimension. It is due to the orientability of real projective
spaces of odd dimension, which implies that the group $PGL(n+
2,\R)$ of real projective transformations of $P^{n+1}$ has two
connected components if $n$ is even. In our case, $n= 4$, so some
of the coarse deformation classes of real projective nonsingular
cubic fourfolds may apriori consist of two deformation classes.

This leads us to a study of the following chirality phenomenon. We
say that a real nonsingular cubic $X\subset P^5$ and its coarse
deformation class are {\it chiral} if $X$ and its mirror image
$X'$ (that is the image of $X$ under a reflection in a hyperplane)
belong to different connected components of $\Cal{C}_\R$, and {\it
achiral} if they belong to the same component (that is if $X$ and
$X'$ can be connected by a continuous family of real non-singular
cubics). Clearly, a coarse deformation class consists of two
deformation classes if and only if it is chiral.

In the present paper we reduce the chirality problem to a specific
problem of the arithmetics of lattices and use this reduction to
show that certain real cubic fourfolds are chiral, while certain
other real cubic fourfolds are achiral. We pay a special attention
to real cubic fourfolds with extremal values of the sum of the
Betti numbers. Namely, we consider in details the cases of {\it
$M$-cubics}, in which $\dim H_*(X(\R);\Bbb Z/2)=\dim
H_*(X(\C);\Bbb Z/2)$ (the maximal value), and the cases of {\it
$(M-1)$-cubics}, in which $\dim H_*(X(\R);\Bbb Z/2)=\dim
H_*(X(\C);\Bbb Z/2)-2$ (the next value). As is shown in \cite{FK},
the $M$-cubics form three and the $(M-1)$-cubics form six coarse
projective classes. In the present paper we prove that one coarse
class of $M$-cubics and three coarse classes of $(M-1)$-cubics are
achiral, while the other coarse classes of $M$- and $(M-1)$-cubics
are chiral.

As a by-product, we give a new proof (in a sense, more natural and
more direct) of the \emph{homological quasi-simplicity} of cubic
fourfolds, where the latter means that two real nonsingular cubic
hypersurfaces $X_1,X_2$ in $P^5$ are coarse deformation equivalent
iff the involutions induced by the complex conjugation on
$H_4(X_i(\C)), i=1,2$, regarded as a lattice via the intersection
index form, are isomorphic (cf. Theorem 1.1 in \cite{FK} and
Proposition 4.1.2 below).

In our previous paper \cite{FK}, we were using a relation between
the nodal cubics in $P^{5}$ and the complete intersections of
bi-degree $(2,3)$ in $P^{4}$. Since these complete intersections
are the $6$-polarized K3-surfaces, it had allowed us to apply
Nikulin's coarse deformation classification of real $6$-polarized
K3-surfaces in terms of involutions on the K3-lattice and his
results on the arithmetics of such involutions, see
\cite{Nikulin},\cite{N2}.

Such a roundabout approach was imposed by a lack of sufficiently
complete understanding of the moduli of cubic hypersurfaces,
contrary to that of K3-surfaces. In particular, in the case of
K3-surfaces one had in one's hands the surjectivity of the period
map, while for cubic fourfolds the characterization of the image
of the period map remained unknown. The situation has changed
recently, after R.~Laza \cite{Laza} and E.~Looijenga
\cite{Looijenga} established a suitable surjectivity statement for
cubic fourfolds.

In our opinion, the two approaches are complementary and both
deserve to be developed further. Combined together they should
give us a better understanding of the topology of the moduli space
of real cubic fourfolds on one hand, and of the topology of the
discriminant of cubic fourfolds on the other hand. Note that
already in \cite{FK} not only the coarse deformation classes but
also their adjacencies were found. Now, via the period map, the
deformation classes become endowed with a certain polyhedral
structure. This opens a way for a full understanding of some
natural stratifications of the moduli and the coefficient spaces
of real cubic fourfolds.

Topological study of nonsingular real cubic hypersurfaces has a
long history, see \cite{FK} for a brief account. In addition, we
would add a reference to the recent investigation of the moduli
space of real cubic surfaces performed by D.~Allcock, J.~Carlson,
and D.~Toledo \cite{ACT}.

Let us recall that according to Klein's classification of real
cubic surfaces, see \cite{Kl} (the classification statement is
reproduced in \cite{FK}), all the real nonsingular cubic surfaces
are achiral. It may be worth mentioning that Klein's achirality
argument in \cite{Kl} contained a mistake, which was corrected by
Klein in his Collected Papers, see \cite{Kl2}.

The paper is organized as follows. In Section 2, we review some
properties of the period map for complex cubic fourfolds. In
Section 3, we introduce the real period spaces with the real
period map and derive the properties of the latter from the
corresponding properties of the complex period map. The results of
Section 3 are applied then in Section 4 to reduce the chirality
problem to some arithmetics of hyperbolic integer lattices and
their reflection groups. Section 5 collects necessary information
about Vinberg's algorithm for finding the fundamental domains of
the arithmetical reflection groups. The technique developed in
Sections 3 -- 5 is applied in Sections 6 and 7 to treat the
chirality of $M$- and $(M-1)$-cubic fourfolds. Section 8 is
devoted to concluding remarks. We mention some other cases which
were studied using similar methods, and mention some other related
results and possible directions of their development. In
particular we discuss a notion of reversibility, which is closely
related to chirality.

\subheading{Acknowledgements} This work has received the first
impulse during the stay of the second author at the Max Planck
Institute in Bonn in August 2007, it took the shape then during
the visit of the first author to the Strasbourg University in
October 2007, and the research was completed during the stay of
the both authors at the Bernoulli Center (EPFL) in Lausanne in the
spring of 2008. We thank all these institution for hospitality and
excellent working conditions. The second author acknowledges a
support from the ANR-05-0053-01 grant of Agence Nationale de la
Recherche.

\heading \S2. Period map for complex cubic fourfolds \endheading

\subheading{2.1. The period domain for marked cubic fourfolds}
Consider a non-singular cubic fourfold $X\subset P^5$. It is well
known that its non-zero Hodge numbers in dimension four are
$h^{3,1}=h^{1,3}=1$ and $h^{2,2}=21$. The lattice $\M(X)=H^4(X)$
is odd with signature $(21,2)$. The {\it polarization class}
$h(X)\in \M(X)$, that is the square of the hyperplane section, is
a characteristic element of $\M(X)$ with $h^2=3$, and so the
primitive sublattice $\M_0(X)=\{x\in \M(X)\,|\,xh=0\}$ is even and
has discriminant group $\Z_3$. This implies that there is a
lattice isomorphism between $\M(X)$ and $\M=3I+2U+2E_8$, which
sends $h(X)$ to $h=(1,1,1)\in 3I$, so that $\M_0(X)$ is identified
with $\M_0=A_2+2U+2E_8$. A particular choice of such an
isomorphism $\phi: (\M(X), h(X))\to (\M,h)$ will be called a {\it
marking of $X$}. We restrict the choice of markings as is
specified below.

The complex line $\phi(H^{3,1}(X))\subset \M_0\otimes\C$ is
isotropic and has negative pairing with the conjugate (and thus,
also isotropic) line
$\phi(H^{1,3}(X))=\overline{\phi(H^{3,1}(X))}$, that is to say,
$w^2=0$, and $w\overline w<0$, (and thus $\overline w^2=0$) for
all $w\in \phi(H^{3,1}(X))$. Writing $w=u+iv$, $u,v\in
\M_0\otimes\R$, we can reformulate it as $u^2=v^2<0$ and $uv=0$,
which implies that the real plane $\la u,v\ra\subset
\M_0\otimes\R$ spanned by $u$ and $v$ is negative definite and
bears a natural orientation given by $u=\Re w, v=\Im w$. Note that
the orientation determined similarly by the complex line
$\phi(H^{1,3}(X))\subset \M_0\otimes\C$ is the opposite one.

The line $\phi(H^{3,1}(X))\subset \M_0\otimes\C$ specifies a point
$\Omega(X)\in P(\M_0\otimes\C)$ (as usual, $P$ states for the
projectivization) called the {\it period point of $X$}. This
period point belongs to the quadric $Q= \{w^2=0\}\subset\PMC$, and
more precisely, to its open subset, $\widehat{\Cal D}=\{w\in
Q\,|\,w\overline w<0\}$.
 This subset has two connected components,
which are exchanged by the complex conjugation (this reflects also
switching from the given complex structure on $X$ to the complex
conjugate one).

The orthogonal projection of a negative definite real plane in
$\M_0\otimes\R$ to another one is non-degenerate. Thus, to select
one of the two connected components of $\widehat{\Cal D}$ we fix
an orientation of negative definite real planes in $\M_0\otimes\R$
so that the orthogonal projection preserves it. We call it the
{\it prescribed orientation} and restrict the choice of markings
to those for which the orientation of $\phi(H^{1,3}(X))$ defined
by the pairs $u=\Re w, v=\Im w$ for $w\in\phi(H^{1,3}(X))$ is the
prescribed one. We denote this component by $\Cal D$ and call it
the {\it period domain}. By $\Aut^+(\M_0)$ we denote the group of
those automorphisms of $\M_0$ which preserve the prescribed
orientation (and thus preserve $\Cal D$). We put
$\Aut^-(\M_0)=\Aut(\M_0)\setminus\Aut^+(\M_0)$. This complementary
coset consists of automorphisms exchanging the connected
components of $\widehat{\Cal D}$.

The projective space $P_{4,3}$ formed by all cubic fourfolds
splits into the {\it discriminant hypersurface} $\D_{4,3}$ formed
by singular cubics and its complement, $\Cal C$. Let $\Cal
C^\sharp$ denote the space of marked non-singular cubics. The
natural projection $\Cal C^\sharp\to\Cal C$ is obviously a
covering with the deck transformation group $\Aut^+(\M_0)$. The
above conventions define the {\it period map} $\per\: \Cal
C^\sharp\to \Cal D$,  $(X,\phi)\mapsto \phi(H^{1,3}(X)$.

\subheading{2.2. Principal properties of the period map} The
global Torelli theorem for cubic fourfolds proved in \cite{Voisin}
claims injectivity of the period map. We need the following
version of this theorem.

\tm{2.2.1. Global Torelli Theorem \cite{Voisin}} Assume that
$(X,\phi)$ and $(X',\phi')$ are non-singular marked cubic
fourfolds such that $\per(X,\phi)=\per(X',\phi')$. Then there
exists one and only one isomorphism $f\:X'\to X$ such that
$\phi'\circ f^*=\phi$.\qed
\endtm

The existence statement is explicit in \cite{Voisin}. The
uniqueness statement is implicit there. It follows easily from two
well known observations: first, each automorphism of a nonsingular
cubic fourfold is induced by a projective transformation, and,
second, if a projective transformation acts trivially in the
cohomology then it is trivial.

\tm{2.2.2. Construction of (anti-)isomorphisms} Let $X$ and $X'$
be non-singular cubic fourfolds and $F: H^4(X;\Z) \to H^4(X';\Z)$
an isometry such that $F(h)=h'$. \roster \item If
$F(H^{3,1}(X))=H^{3,1}(X')$, then there exists one and only one
isomorphism $f\:X'\to X$ such that $f^*=F$. \item If
$F(H^{1,3}(X))=H^{3,1}(X')$, then there exists one and only one
isomorphism $f\:X'\to\overline X$ such that $f^*=F$.
\endroster
\endtm

Here and in what follows we denote by $\overline X$ the variety
complex conjugate to $X$. If $X\subset P^5$ is given by a
polynomial, then $\overline X\subset P^5$ can be seen as the
variety given by the polynomial with the complex conjugate
coefficients.

\demo{Proof of 2.2.2} The first statement is nothing but an
equivalent version of Theorem 2.2.1. The second statement follows
from the first one or directly from Theorem 2.2.1 applied to
$(X',\phi')$ and $(\overline X,\phi'\circ F)$, where $\phi'$ is
any marking of $X'$.\qed
\enddemo

Consider the reflection $R_v$ in $\MC$ across the
mirror-hyperplane $H_v=\{x\in \MC\, |\, xv=0\}$ defined as
$x\mapsto x-2\frac{xv}{v^2}v$, and note that it preserves the
lattice $\M_0$ invariant if $v\in \M_0$ is such that $v^2=2$, or
such that $v^2=6$ and $xv$ is divisible by $3$ for all $x\in
\M_0$. We call these two types of lattice elements {\it $2$-roots}
and {\it $6$-roots} respectively, and denote their sets by $V_2$
and $V_6$. Note that $R_v\in\Aut^+(\M_0)$ for any $v\in V_2\cup
V_6$. If $v\in V_2$, then the reflection $R_v$ extends (as a
reflection) to $\M$ and $h$ is preserved by this extension. By
contrary if $v\in V_6$, the reflection $R_v$ does not extend to a
reflection in $\M$, and moreover, the unique extension of $R_v$ to
$\M$ maps $h$ to $-h$ (cf. Lemma 4.3.2 below). On the other hand,
if $v\in V_6$ then the \emph{anti-reflection} $-R_v$ extends to an
isometry of $\M$ preserving $h$. This extension is the
anti-reflection with respect to the 2-plane generated by $h$ and
$v$. In particular, it represent also an element of
$\Aut^+(\M_0)$.

The union of the mirrors $H_v$ for all $v\in V_2$ gives after
projectivization a union of hyperplanes $ \Cal H_\D\subset\PMC$,
and a similar union for all $v\in V_6$ gives another union of
hyperplanes, $\Cal H_\infty \subset\PMC$.

\tm{2.2.3. Surjectivity of the period map
\cite{Looijenga},\cite{Laza}} The image of the period map $\per\:
\Cal C^\sharp\to \Cal D$ is the complement of $\Cal H_\D\cup\Cal
H_\infty$.\qed
\endtm

According to the Griffiths theory, for any nonsingular cubic
$X\subset P^5$ the line $H^{3,1}(X)$ is spanned by the class
$[\omega_p]\in H^4(X;\C)$ of the 4-form $\omega_p=\Res(\Cal
E/p^2)$. Here $\Cal E$ stands for the Euler $5$-form in $\C^6$,
$\Cal E=\sum_{i=0}^5(-1)^i x_i \d x_0\wedge \dots \wedge \d
\widehat x_i \wedge \dots\wedge \d x_5$, and $p$ for a polynomial
defining $X$ (as usual, a hat over $x_i$ means that this term is
omitted). The ratio $\Cal E/p^2$ is a well-defined meromorphic
5-form in $P^5$, with a second order pole along $X$. The residue
$\omega_p$ of this form is a $4$-form on $X$, which is a linear
combination of $(3,1)$ and $(4,0)$-forms. Its class $[\omega_p]$
is known to be non-trivial, thus, it spans $H^{3,1}(X)$.

\heading \S3. Periods in the real setting \endheading

\subheading{3.1. Geometric involutions} Consider a non-singular
cubic fourfold $X$ defined by a real polynomial $p$, and let
$\conj\:X\to X$ denote the complex conjugation map, which will be
called also the {\it real structure on $X$}. The latter map
induces a lattice involution $\conj^*\:\M(X)\to \M(X)$ such that
$\conj^*(h)=h$ and, hence, induces also a lattice involution in
$\M_0(X)$. Denote by $\M^\pm_0(X)$ and $\M^\pm(X)$ the
eigen-sublattices $\{x\in \M_0(X)\, |\, \conj^*(x)=\pm x\}$ and
$\{x\in \M(X)\, |\, \conj^*(x)=\pm x\}$, respectively. We have
obviously $\M^-=\M^-_0$ and
$\sigma_-(\M^+(X))=\sigma_-(\M^+_0(X))$, where  $\sigma_-$ denotes
the negative index of inertia.

\tm{3.1.1. Lemma} One has $\sigma_-(\M^\pm_0(X))=1$.
\endtm

\demo{Proof} The map $w\mapsto\overline{\conj^*w}$ gives an
anti-linear involution in $H^{3,1}(X)$. Thus, there exist non zero
elements $w\in H^{3,1}(X)$ such that $\conj^*(w)=\overline w$. In
terms of the real and imaginary components of $w=u_++iu_-$, this
identity means that $u_\pm\in \M^\pm(X)\otimes\R$. These
components satisfy the relations $u_+^2=u_-^2=\frac12 w\overline
w<0$. They belong to $\M_0^\pm(X)$, since $wh=0$. It remains to
notice that the intersection form is positive definite on
$H^{2,2}(X)$.\qed\enddemo

We call a lattice involution $c:\M\to \M$ {\it geometric } if
$c(h)=h$ and $\sigma_-(\M^\pm_0(c))=1$, where $\M^\pm_0(c)$
denotes the eigen-sublattices $\{x\in \M_0 \,|\,c(x)=\pm x\}$. Let
us note that all geometric involutions preserve $\M_0$ and the
involutions induced in $\M_0$ belong to $\Aut^-(\M_0)$.

According to Lemma 3.1.1, all lattice involutions $c:\M\to \M$
isomorphic to an involution $\conj^*\:\M(X)\to \M(X)$ for a
non-singular real cubic $X$ are geometric. A pair $(c:\M\to \M,
h\in \M)$ isomorphic to  $(\conj^*\:\M(X)\to \M(X), h(X))$ is
called the {\it homological type} of $X$. By a {\it real
$c$-marked cubic fourfold} we understand a real non-singular cubic
fourfold equipped with a marking $\phi$ such that
$\phi\circ\conj^*=c\circ\phi$.

\tm{3.1.2. Theorem} For any geometric involution $c:\M\to \M$ the
pair $(c,h)$ is the homological type of some non-singular real
cubic fourfold.
\endtm

This theorem is one of the results obtained in \cite{FK}. After
fixing some notation, we will give below (at the end of subsection
3.3) an independent proof based on the surjectivity of the period
map and the global Torelli theorem.

The number of isometry classes of geometric involutions is finite.
Their list can be found in \cite{FK} (see also tables 8,9 in \S8).

Up to the end of this section we suppose that $c$ is a geometric
involution.

\subheading{3.2. Real parameter space $\Cal C_\R^{c}$} We denote
by $\Cal C^c_\R\subset\Cal C_\R$ the set of real cubic fourfolds
of homological type $c$, and by $\Cal C_\R^{c\sharp}$ the set of
$c$-marked real cubic fourfolds. The former consists of some
number of connected components of $\Cal C_\R$. The latter can be
seen as the real part of $\Cal C^\sharp$ with respect to the
involution $\conj^{c\sharp}\:\Cal C^\sharp\to\Cal C^\sharp$, which
send $(X,\phi)\in\Cal C^\sharp$ to
$(\conj(X),c\circ\phi\circ\conj^*)$. The forgetful map
$(X,\phi)\to X$ defines a (multi-component) covering $\Cal
C_\R^{c\sharp}\to\Cal C_\R^{c}$ with the deck transformation group
$\Aut^+(\M_0)$.

\subheading{3.3. Real period domain $\Cal D^c_\R$} Let us extend
$c$ to a complex linear involution on $\M\otimes \C$ and denote
also by $c$ the induced involutions on $\MC$, $P=P(\MC)$, and
$\widehat{\Cal D}$. Note that $c$ permutes the two components
$\Cal D$ and $\overline{\Cal D}$ of $\widehat{\Cal D}$, and thus,
$\overline c(\Cal D)=\Cal D$, where $\overline c\:\MC\to\MC$ is
the composition of $c$ with the complex conjugation in $\MC$.

Let $\widehat{\Cal D}_\R^{c}$ and $\Cal D_\R^{c}$ denote the fixed
point set of $\overline c$ restricted to $\widehat{\Cal D}$ and
$\Cal D$. The latter consists of the lines generated by
$w=u_++iu_-$ such that $u_\pm\in \M_0^\pm(c)\otimes\R$,
$u_+^2=u_-^2<0$, and the orientation $u_+,u_-$ is the prescribed
one. Since $c$ is geometric, both ${\Cal D}_\R^{c}$ and its
(trivial) double covering $\widehat{\Cal D_\R^{c}}$ are nonempty.

As it follows from definitions, the period of a c-marked real
cubic fourfold belongs to $\Cal D_\R^{c}=\{x\in\Cal
D\,|\,c(x)=\overline x\}$. Therefore, we call $\Cal D_\R^{c}$ the
{\it real period domain of real $c$-marked cubic fourfolds}.

\demo{Proof of Theorem 3.1.2} Pick up a generic point $[w]\in
{\Cal D}_\R^c $ (so that there is no vector $v\in V_2\cup V_6$
orthogonal to $w$) and apply Theorem 2.2.3. This gives a
nonsingular cubic fourfold $X$ and a marking $\phi$ such that
$\per(X,\phi)=[w]$. The triple $(X,X'=\overline X,
F=\phi^{-1}c\phi)$ satisfies the assumptions of Theorem 2.2.2,
which gives an antiholomorphic involution $\conj:X\to X$ such that
$\conj^*= \phi^{-1}c\phi$. Clearly, $(\M,c)$ is the homological
type of $(X,\conj)$, and it remains to notice that $\Pic X=\Z$,
$X(\R)$ is non empty (as it is for any real hypersurface of odd
degree), and therefore any antiholomorphic involution of $X$ is
induced by the complex conjugation in suitable projective
coordinates of $P^5=P(\Cal O_X(1))$.\qed
\enddemo

\subheading{3.4. Refined real period map} Consider the quadratic
cones $\Upsilon_\pm(c)=\{u\in \M_0^\pm(c)\otimes\R \: u^2<0\}$ and
the associated Lobachevsky (one- and two-component, respectively)
spaces $\Lambda_\pm(c)=\Upsilon_\pm(c)/\R^*$ and
${\Lambda}^\sharp_\pm(c)=\Upsilon_\pm(c)/\R_+$, where
$\R^*=\R\smallsetminus\{0\}$ and $\R_+=(0,\infty)$.

Like in 3.3, we associate with a point in $\Cal D_\R^c$
represented by $w=u_++iu_-$ {\rm (where $u_\pm\in
\M_0^\pm(c)\otimes\R$, $u_+^2=u_-^2<0$, and the orientation
$u_+,u_-$ is the prescribed one)} the point in
$\Lambda_+(c)\times\Lambda_-(c)$ represented by the pair
$(u_+,u_-)$. This gives a well-defined analytic isomorphism $\Cal
D_\R^c=\Lambda_+(c)\times\Lambda_-(c)$. The ambiguity in the
choice of representatives gives rise to a refined real period
domain $\Cal D_\R^{c\sharp}\subset
{\Lambda}^\sharp_+(c)\times{\Lambda}^\sharp_-(c)$, $\Cal
D_\R^{c\sharp}=\{(u_+\R_+, u_-\R_+)\in
{\Lambda}^\sharp_+(c)\times{\Lambda}^\sharp_-(c)\,\vert\,$ the
orientation $ u_+,u_- $ is the prescribed one$\}.$

To define $ {\per}_\R^{c\sharp}(X,\phi)\in\Cal D_\R^{c\sharp}$ for
a non-singular real $c$-marked cubic $(X,\phi)\in\Cal
C^{c\sharp}_\R$, we pick up a real polynomial $p$ defining $X$ and
consider $w=\phi([\omega_p])$ (see the end of Section 2). As we
have seen already, the latter splits as $w=u_++iu_-$, where
$u_\pm\in \M_0^\pm(c)$, the pair $(u_+,u_-)$ is defined uniquely
by $X$ up to a positive factor, and this pair spans a negative
definite plane with the prescribed orientation. Thus, we obtain a
uniquely defined real period ${\per}_\R^{c\sharp}(X,\phi)\in\Cal
D_\R^{c\sharp}$ and a well defined map ${\per}_\R^{c\sharp}\:\Cal
C_\R^{c\sharp}\to\Cal D_\R^{c\sharp}$. The above components
$u_\pm\R_+\in\Lambda^\sharp_\pm(c)$ of
${\per}_\R^{c\sharp}(X,\phi)$ will be denoted
$u_\pm^\sharp(X,\phi)$.

\subheading{3.5. Polyhedral period cells} Denote by $\Cal
H_\pm(c)\subset \Lambda_\pm(c)$ and $\Cal H_\pm^\sharp(c)\subset
\Lambda_\pm^\sharp(c)$ the union of hyperplanes orthogonal to
vectors from $(V_2\cup V_6)\cap \M^\pm_0(c)$. The connected
components of the complement $\Lambda_\pm(c)\setminus \Cal
H_\pm(c)$ will be called {\it the cells of $\Lambda_\pm(c)$} and
the hyperplanes from $\Cal H_\pm(c)$ {\it the walls}. As is known,
these hyperplanes form a locally finite arrangement (the group
generated by reflections in these hyperplanes is discrete) so that
the above cells are (locally finite) polyhedra. Put
$$
{\Per}_\R^{c}=\Cal
D_\R^{c\sharp}\cap((\Lambda_+^\sharp(c)\setminus \Cal
H_+^\sharp(c))\times(\Lambda_-^\sharp(c)\setminus \Cal
H_-^\sharp(c))).
$$
and call {\it $c$-cells} the connected components of
${\Per}_\R^{c}$. Note that the orientation restriction involved in
the definition of $\Cal D_\R^{c\sharp}$ establishes a one-to-one
correspondence between the halves of $\Lambda_+^\sharp(c)$ and the
halves of $\Lambda_-^\sharp(c)$, and this correspondence commutes
with multiplication by $-1$. Therefore, ${\Per}_\R^{c}$ spits into
a union of pairs of opposite $c$-cells. The natural projection
${\Per}_\R^{c}\to \Lambda_+(c)\times\Lambda_-(c)$ establishes a
one-to-one correspondence between the set of pairs of opposite
$c$-cells and the set of products of the cells of
$\Lambda_\pm(c)$.

Given a continuous family of real $c$-marked cubics
$(X_t,\phi_t)$, $t\in [0,1]$, they can be defined by a continuous
family of polynomials $p_t$, and hence their real periods
$u_\pm^\sharp(X_t,\phi_t)$ belong to the same cells of
$\Lambda_\pm^\sharp(c) $. The converse is also true.

\tm{3.5.1. Lemma} Assume that $(X_i,\phi_i)$, $i=0,1$ is a pair of
real $c$-marked cubic fourfolds defined by real polynomials $p_i$.
Then, $X_i$ can be connected by a continuous family $X_t$ of real
$c$-marked cubic fourfolds if and only if their periods
$u_\pm^\sharp(X_i,\phi_i)$ belong to the same cells of
$\Lambda_\pm^\sharp(c)$ {\rm (or in other words, if and only if
the periods $\per_\R(X_i,\phi_i)$ belong to the same component of
$\Per_\R^{c}$).}
\endtm

\demo{Proof} It follows from the description of the periods of
cubic fourfolds (and the local Torelli theorem over the reals),
because the vectors $v\in (V_2\cup V_6)$ which are not from
$\M^+_0\cup \M^-_0$ define hyperplanes $H_v$ which have
intersection with $\M^\pm_0\otimes\R$ of codimension less than
one.\qed
\enddemo

\heading \S4. Deformations and chirality\endheading

\subheading{4.1. The mirror pairs of markings} Given a real
hypersurface $X\subset P^5$, we can consider its {\it mirror
image}, $X'=R(X)$, obtained from $X$ by a reflection $R\:P^5\to
P^5$ with respect to some real hyperplane $H\subset P^n$.
According to our definitions, $X$ is {\it chiral} if $X$ and $X'$
belong to different connected components of $\Cal C_\R$, and {\it
achiral} if they belong to the same component.

Assume that $(X,\phi)$ is a marked non-singular cubic fourfold.
Then the isomorphism $R^*\:\M(X')\to \M(X)$ induced by $R$
respects the Hodge structure and the polarization classes of $X$
and $X'$, and thus yields a marking $\phi\circ R^*$ of $X'$. We
say that the markings $\phi$ and $\phi'=\phi\circ R^*$ are {\it
mirror images of each other}, or {\it a mirror pair of markings}.

\tm{4.1.1. Lemma} Assume that a nonsingular real cubic fourfold
$X$ is defined by a real polynomial $p$ and its mirror image,
$X'$, by a polynomial $q$. Then the period vectors
$\phi[\omega_p]$ and $\phi'[\omega_{q}]$ are oppositely directed
if $X$ and $X'$ are endowed with the mirror pair of markings:
$\phi$ and $\phi'=\phi\circ R^*$.
\endtm

\demo{Proof} The form $\Cal E/q^2$ representing $[\omega_q]$
changes the direction under the action of $R$, because $R^*(\Cal
E)=-\Cal E$ and $q\circ R$ differs from $p$ by a real factor.\qed
\enddemo

As an immediate corollary of Lemma 4.1.1 and Lemma 3.5.1 we get a
new proof of the following theorem from \cite{FK}.

\tm{4.1.2. Coarse deformation classification} One real
non-singular cubic fourfold is deformation equivalent to a
projective transformation of another real non-singular cubic
fourfold if and only if they are of the same homological type.
\endtm

\demo{Proof} Given a $c$-marking, we can compose it with lattice
reflections $R_v, v\in V_2\cap \M_\pm^0(c)$, and anti-reflections
$-R_v, v\in V_6\cap \M_\pm^0(c)$, to move the period into any pair
of opposite cells of $\Per_\R(c)$ given in advance. When
necessary, we can apply  Lemma 4.1.1 and move the period into any
of these opposite cells. According to Lemma 3.5.1 it means that
the real non-singular cubics of homological type $c$ are coarse
deformation equivalent to each other. The "only if" part is
trivial. \qed
\enddemo

\subheading{4.2. Basic criterion of chirality for cubic fourfolds}
Let us fix a geometric involution $c$. Given a nonsingular
$c$-marked real cubic fourfold $(X,\phi)$, denote by
$P^\sharp(X)\subset{\Per}_\R^{c}$ the $c$-cell which contains
${\per}_\R^{c\sharp}(X,\phi)$ (in other words, the $c$-cell which
contains $w=\phi[\omega_p]$ where, as usual, $p$ is a real
polynomial defining $X$).

\tm{4.2.1. Lemma} The underlying nonsingular real cubic fourfold
$X$ of a real $c$-marked cubic fourfold $(X,\phi)$ is achiral if
and only if there exists a lattice isometry of $\M$ which (1)
commutes with $c$, (2) preserves the polarization class $h$, (3)
induces an automorphism of $\M_0$
 which preserves the prescribed orientation,
and (4) sends the $c$-cell $P^\sharp(X)$ to the opposite $c$-cell,
$-P^\sharp(X)$.
\endtm

\demo{Proof} Let $X'$ denote the mirror image of $X$ with the
mirror image marking $\phi'$. By Lemma 4.1.1, its period
$w'=\phi'[\omega_q]$ belongs to $-P^\sharp(X)$. On the other hand,
any continuous family of real nonsingular cubic fourfolds
connecting $X$ with $X'$ gives another marking of $X'$, say
$\phi''$, and according to Lemma 3.5.1 the period
$\phi''[\omega_q]$ belongs to $P^\sharp(X)$. Comparing the two
markings of $X'$ we obtain a lattice isometry of $\M=\M(X')$ which
transforms $P^\sharp(X)$ into $-P^\sharp(X)$; being a difference
between two markings, it also preserves the polarization $h$,
induces an automorphism of $\M_0$ which preserves the prescribed
orientation, and commutes with $c$. Conversely, given such a
lattice isometry, we can change the mirror image marking of $X'$
and then apply Lemma 3.5.1 to deduce that $X$ and $X'$ belong both
to the same component of $\Cal C_\R$.\qed
\enddemo

\subheading{4.3. Few lattice gluing lemmas} To simplify the above
criterion and to reduce it to a study of $\Aut \M_+^0(c)$ we need
the following results involving a technique of discriminant
groups. Recall that for any non-degenerate lattice $\LL$ of finite
rank the discriminant group $\discr \LL=\LL/\LL^*$ is a finite
group and that, if the lattice $\LL$ is even, this group carries a
canonical finite quadratic form $\frak q_\LL:\discr \LL\to\Q/2\Z$
defined via $\frak q_\LL(x+\LL)=x^2\mod 2\Z$. Note that any
isometry, $f\in\Aut \LL$, induces an automorphism of $\discr \LL$,
which preserves $\frak q_\LL$ if $\LL$ is even. This induced
automorphism will be denoted by $\delta(f)$.

\tm{4.3.1. Nikulin's theorem \cite{Nikulin}} Assume that $\LL$ is
an even lattice of signature $(n,1)$, $n\ge0$, whose discriminant
group $\discr(\LL)$ is $2$-periodic. Then any isometry
$\delta\:\discr(\LL)\to\discr(\LL)$ is induced by some isometry
$f\:\LL\to \LL$.\qed
\endtm

In the present paper we deal with the three lattices: $\M_-(c)$,
$\M^0_+(c)$, and the rank $1$ lattice $\langle h\rangle\subset \M$
generated by $h$. The first two lattices are even, and the latter
is odd. The discriminant group $\discr \M_-(c)$ is 2-periodic, the
discriminant group $\discr \langle h\rangle$ is a cyclic group of
order $3$, and the discriminant group $\discr \M^0_+(c)$ is
canonically isomorphic to the direct sum $\discr \M_-(c)+\discr
\langle h\rangle$, so that $\discr \M_-(c)$ is identified with the
2-primary part $\discr_2 \M^0_+(c)$ of $\discr \M^0_+(c)$, and
$\discr \langle h\rangle$ with its 3-primary part $\discr_3
\M^0_+(c)$. The canonical isomorphism $\discr_2 \M^0_+(c)\to
\discr \M_-(c)$ is an anti-isometry, that is it transforms $-\frak
q_{\M_-(c)}$ into $\frak q_{\M^0_+(c)}$ restricted to $\discr_2
\M^0_+(c)$.
 (In fact, the lattice
$\discr\langle h\rangle$, as any non-degenerate finite rank
lattice with a fixed characteristic element, can be also equipped
with a quadratic form, and with respect to this quadratic form the
canonical isomorphism $\discr_3 \M^0_+(c)\to \discr\langle
h\rangle$ is also an anti-isometry.)

The following lattice gluing lemmas are well known and their
proofs are straightforward, see, e.g., \cite{Nikulin}.

\tm{4.3.2. Lemma} Any automorphism $f_+^0\in\Aut(\M_+^0(c))$ can
be uniquely extended to $\M_+(c)$. This extension sends the
polarization class $h$ to itself if and only if the 3-primary
component $\delta_3(f_+^0)$ of $\delta(f_+^0)$ is trivial, that is
$\delta_3(f_+^0)=\id$.\qed\endtm

\tm{4.3.3. Lemma} A pair of automorphisms
$f_\pm\in\Aut(\M_\pm(c))$ are induced from $f\in\Aut(\M,c)$ if and
only if $\delta(f_+)=\delta(f_-)$.\qed
\endtm

Automorphisms $f_\pm$ satisfying the conditions of Lemma 4.3.3
will be called {\it compatible}.

\subheading{4.4. Lattice characterization of chirality} The
\emph{reflection group} $W_+$ generated in $\Aut(\M_+^0(c))$ by
reflections $R_v$, $v\in (V_2\cup V_6)\cap \M_+^0(c)$ acts
transitively on the set of cells of $\L_+(c)$. If $v\in V_6$, then
$R_v$ does not extends to $M_+$, but anti-reflection $-R_v$ does.
So, we consider also the group $W_+^\# \subset$ $\Aut(\M_+^0(c))$
generated by reflections $R_v$, $v\in V_2\cap \M_+^0(c)$, and
anti-reflections $-R_v$, $v\in V_6\cap \M_+^0(c)$ (the two groups
are isomorphic and induce the same action on $\L_+$). Any of the
cells $P_+\subset\L_+(c)$ being fixed, the group $\Aut(\M_+^0(c))$
splits into a semi-direct product $W_+\rtimes\Aut(P_+)$, where
$\Aut(P_+)=\{g\in\Aut(\M_+^0(c))\,|\,g(P_+)=P_+\}$ is the
stabilizer of $P_+$.

With $\M_-(c)$ the situation is even simpler: since its
discriminant group is of period $2$ the intersection $V_6\cap
\M_-(c)$ is empty. Thus, in this case we consider simply the
reflection group $W_-\subset\Aut(\M_-(c))$ generated by
reflections $R_v$, $v\in V_2\cap \M_-^0(c)$. This reflection group
acts transitively on the set of cells of $\L_-(c)$ and, therefore,
$\Aut(\M_-(c))$ splits into a semi-direct product
$W_-\rtimes\Aut(P_-)$, where $\Aut(P_-)=\{g\in\Aut(\M_-(c))|
g(P_-)=P_-\}$ is the stabilizer of a cell $P_-$ of $\L_-(c)$.

The preimage of $P_\pm$ in $\L^\#_\pm$ is the union of a pair of
$c$-cells, $P_\pm^\#$ and $-P_\pm^\#$. Each $g\in\Aut(P_\pm)$
either permutes this pair of cells, and then we say that it is
{\it $P_\pm$-reversing}, or it preserves both $P_\pm^\#$ and
$-P_\pm^\#$, and then we call it {\it $P_\pm$-direct}. The
subgroup of $\Aut(P_\pm)$ formed by $P_\pm$-direct elements will
be denoted by $\Aut^+(P_\pm)$, while the coset of
$P_\pm$-reversing elements will be denoted by $\Aut^-(P_\pm)$. The
crucial for our study of chirality observation is that an
automorphism $f\in\Aut(\M)$ preserving each of $P_\pm$ belongs to
$\Aut^+(\M)$ if and only if its components $f_+=f|_{\M_+^0}$,
$f_-=f|_{\M_-}$ are both of the same type:
 either simultaneously $f_\pm\in\Aut^+(P_\pm)$,
or simultaneously $f_\pm\in\Aut^-(P_\pm)$.

In the case of lattices $\M_+^0$, an additional characteristic of
$g\in\Aut(\M_+^0)$ is its 3-primary component, $\delta_3(g)$,
which may be trivial or not. In a bit more general setting, we
consider a hyperbolic lattice $\LL$ whose discriminant splits as
$\discr(\LL)=\discr_2(\LL)+\discr_3(\LL)$, where $\discr_2(\LL)$
is 2-periodic and $\discr_3(\LL)=\Z/3$. We say that
$g\in\Aut(\LL)$ is {\it $\Z/3$-direct} if $\delta_3(g)=\id$, and
{\it $\Z/3$-reversing} if $\delta_3(g)\ne\id$ (certainly, in the
later case $\delta_3(g)=-\id$).

 \tm{4.4.1. Theorem} A non-singular real cubic fourfold
$X$ of homological type $c$ is achiral if and only if the lattice
$\M^0_+(c)$ admits an automorphism $g\in\Aut^-(P_+)$ which is
$\Z/3$-direct.
\endtm

\demo{Proof} The "only if" part is a straightforward consequence
of the "only if" part of Lemma 4.2.1.

To prove the "if" part, let us pick up a $c$-marking
$\phi\:\M(X)\to \M$ and choose $f_+^0\in\Aut^-(P_+(X))$ which is
$\Z/3$-direct. From Lemma 4.3.2 it follows that $f_+^0$ extends to
$f_+\in\Aut \M_+$ preserving $h$. Lemma 4.3.3 and Theorem 4.3.1
imply that we can find $f_-\in\Aut(\M_-)$ compatible with $f_+$
and $f\in\Aut(\M)$ defined by $(f_+,f_-)$. By composing $f_-$ (and
$f$) with a suitable $w_-\in W_-$, the component $f_-$ can be
chosen in $\Aut(P_-)\subset\Aut(\M_-)$. If $f\in\Aut(\M)$ defined
by $(f_+,f_-)$ belongs to $\Aut^+(\M)$, then $f$ transforms
$P^\sharp(X)$ into $-P^\sharp(X)$ since it preserves the
prescribed orientation and $f_+\in\Aut^-(P_+)$. Therefore, in this
case due to Lemma 4.2.1 we are done. If $f\in\Aut^-(\M)$, then we
replace $f_-$ by $-f_-$, observe that the pair $(f_+,-f_-)$
defines an automorphism $f\circ c$ which belongs to $\Aut^+(\M)$,
and argue as before. \qed
\enddemo

\heading \S5. Auxiliary arithmetics
\endheading

\subheading{5.1. Root systems and chirality of special hyperbolic
lattices} In this section $\LL$ is a lattice of signature $(n,1)$,
$n\ge1$. Throughout this section we make two additional
assumptions on $\LL$ which are satisfied in the cases of
$\LL=\M_+^0(c)$ that we are concerned about.
 The first assumption is that the discriminant
$\discr(\LL)$ splits as $\Z/3+\discr_2(\LL)$, where the summand
$\discr_2(\LL)$ is a $2$-periodic group. Let $\Phi=V_2\cup V_6$,
where $V_k=\{v\in \LL\,|\,v^2=k,2\,\frac{vw}{v^2}\in\Z,\, \forall
w\in \LL\}$ (note that for $k=2$ the condition $2
\frac{vw}{v^2}\in\Z$ is always satisfied). Our second assumption
is that the rank of $\Phi$ equals to the rank of $\LL$ (that is
maximal possible), and, thus, $\Phi$ is a root system in $\LL$.
This holds for $\LL=\M_+^0(c)$ for all geometric involution $c$
except one rather special case $\M_+^0(c)=U(2)+E_6(2)$ in which
$\Phi=\oo$ (the complete list of $\LL=\M_+^0(c)$ is given in
Tables 8--9, in section 8, see also \cite{FK} for more details).
Vectors $v\in V_k$ will be called {\it $k$-roots}.

We let $\LL_\R=\LL\otimes\R$, and like before, consider
$\Upsilon=\{v\in \LL_\R\,|\,v^2<0\}$, and the hyperbolic spaces
$\L=\Upsilon/\R^*$, along with $\L^\#=\Upsilon/\R_+$. In this
context we use notation $H_v$ for the hyperplane $\{w\in
\LL_\R\,|\,vw=0\}$ and $H_v^\pm$ for the half-spaces $\{w\in
\LL_\R\,|\,\pm vw\ge0\}$. For $v\in\Upsilon$, $H_v$, etc., we
denote by $[v]\in\L$, $[H_v]\subset\L$, $[v]^\#\in\L^\#$,
$[H_v]^\#\subset\L^\#$, etc., the corresponding object after
projectivization.

We distinguish the {\it reflection group} $W\subset\Aut(\LL)$
generated by the reflections $R_v\in\Aut(\LL)$, $x\mapsto
x-2\frac{vx}{v^2}v$, $v\in V_2\cup V_6$, and the group
$W^\#\subset\Aut(\LL)$ generated by the reflections $R_v, v\in
V_2$, and the anti-reflections $-R_v, v\in V_6$. Hyperplanes
$[H_v]$ (respectively $[H_v]^\#$), $v\in\Phi$, cut $\L$
(respectively $\L^\#$) into open polyhedra, whose closures are
called the {\it cells}. The cells in $\L$ are the fundamental
chambers of $W$, and the pairs of opposite cells in $\L^\#$ are
the fundamental chambers of $W^\#$.

Let us pick up a cell $P\subset\L$ and fix a covering $c$-cell
$P^\#\subset\L^\#$. Choosing any vector $p\in\Upsilon$ so that
$[p]^\#$ lies in the interior of $P^\#$, we let
 $\Phi^\pm=\{v\in \Phi|\pm vp>0\}$. The minimal subset $\Phi^b\subset\Phi^-$
such that $P^\#=\cap_{v\in\Phi^b}\,[H_v^-]^\#$ is called the {\it
basis of $\Phi$ defined by $P^\#$}. The hyperplanes $[H_v]$,
$v\in\Phi^b$, support $n$-dimensional faces of $P$ and will be
called the {\it walls of $P$}. Note that any $v\in \Phi^-$ is a
linear combination of the roots in $\Phi^b$ with non-negative
coefficients.

Theorem 4.4.1 motivates the following definition: $\LL$ is called
\emph{achiral} if it admits a $\Z/3$-direct automorphism
$g\in\Aut^-(P)$, for some cell $P$. Obviously, if $\LL$ is achiral
then a $\Z/3$-direct automorphism $g\in \Aut^-(P)$ exists for any
cell $P$. It is also obvious that existence of a $\Z/3$-direct
$g\in\Aut^-(P)$ is equivalent to existence of $\Z/3$-reversing
$h\in\Aut^+(P)$, since these two kinds of automorphisms just
differ by sign.

\subheading{5.2. Coxeter's graphs and their symmetry} The {\it
Coxeter graph} $\G$ has $\Phi^b$ as the vertex set. The vertices
are colored: $2$-roots are white and $6$-roots are black. The
edges are weighted: the weight of an edge connecting vertices
$v,w\in\Phi^b$ is $m_{vw}=4\frac{(vw)^2}{v^2w^2}$, and $m_{vw}=0$
means absence of an edge. These weights are non-negative integers,
because $2\frac{ vw}{v^2},2\frac{vw}{w^2}\in\Z$, and $v^2,w^2>0$
for any $v,w\in\Phi^b$. In the case of $m_{vw}=1$, the angle
between $H_v$ and $H_w$ is $\pi/3$, and $v^2=w^2$; such edges are
not labelled.
 The case of $m_{vw}=2$ (corresponds to angle $\pi/4$) cannot
happen, since $v^2,w^2\in\{2,6\}$. An edge of weight $m_{vw}=3$
connects always a $2$-root with a $6$-root; it corresponds to
angle $\pi/6$, and will be labelled by $6$.
 The case of $m_{vw}=4$ corresponds to parallel hyperplanes in $\L$, and we sketch a
thick edge between $v$ and $w$. If $m_{vw}>4$, then  the
corresponding hyperplanes in $\L$ are ultra-parallel (diverging),
and we sketch a dotted edge.

For a subset $J\subset \Phi^b$ we may consider also the subgraph
$\G_J$ which is formed by the vertex set $J$ and all the edges of
$\G$ connecting these vertices. We say that $\G_J$ is the {\it
Coxeter graph of $J$}. If $J$ is finite and ordered,
$J=\{v_1,\dots,v_{|J|}\}$, then we consider also the {\it Gram
matrix}, $G_J$, whose $(ij)$-entry is $v_iv_j$.

A permutation $\s\:J\to J$ will be called a {\it symmetry of
$\G_J$} if it preserves the weight of edges and the length of the
roots, i.e., $(\s(v))^2=v^2$ and $m_{\s(v)\s(w)}=m_{vw}$ for all
$v,w\in J$.

\tm{5.2.1. Existence of symmetries} Assume that a subset
$J\subset\Phi^b$ spans $\LL$ over $\Z$. Then any symmetry
$\s\:J\to J$ of $\G_J$ is induced by an automorphism of the
lattice $\LL$ which preserves the cell $P^\#$ invariant.
\endtm

\demo{Proof} Such a symmetry preserves the Gram matrix of the
vectors from $J$. Therefore, it is induced by an isometry of
$\LL\otimes\Q$. Since the vectors from $J$ span $\LL$ over $\Z$,
this isometry maps $\LL$ to $\LL$. Assuming that it maps $P^\#$ to
another cell, we observe that these two cells have $J$ as a common
set of face normal vectors. Pick up a wall separating the two
cells and notice that each of the normal root vectors $\pm v\ne 0$
of such a wall has non-negative product with the vectors from $J$,
which is a contradiction, since the vectors from $J$ generate the
whole space.\qed
\enddemo

To recognize $\Z/3$-reversing symmetries of $\G$, one can use the
following observation. Considering some direct sum decomposition
of $\LL$, we observe that one of the direct summands, $\LL_1$, has
$\discr_3(\LL_1)=\Z/3$, while the other direct summands have
$2$-periodic discriminants (because $\discr(\LL)$ gets an induced
direct sum decomposition). For any vertex $w$ of $\G$ viewed as a
vector in $\LL$, we can consider its $\LL_1$-component. Our simple
observation is that {\it $\s$ is $\Z/3$-direct if for all black
vertices, $v\in V_6$, of $\G$ the $\LL_1$-components of $v$ and
$\s(v)$ are congruent modulo $3\LL_1$, and $\Z/3$-reversing if for
some $v\in V_6$ we have $v-\s(v)\notin3\LL$.}

\subheading{5.3. Vinberg's algorithm} Vinberg's method
\cite{Vinberg1} of calculation the Coxeter graph of $\Phi$ is to
pick up a vector $p\in\Upsilon$ so that $[p]^\#\in P^\#$, and then
to determine a sequence of roots $v_i\in\Phi^b$, $i=1,2,\dots$,
ordered so that the hyperbolic distance from $p$ to the walls
$H_i=H_{v_i}$ of $P$ is increasing. Such distance can be
characterized by the (non-negative) value
$d_i=d_i(p)=2\frac{(pv_i)^2}{(v_i)^2}$, which will be called {\it
the level of root $v_i$ with respect to $p$} (the coefficient $2$
here is chosen to make $d_i$ integer in the further
considerations).

The level zero vectors in Vinberg's sequence form a root basis in
the root system $\{v\in V\vert vp=0\}$. Since choosing of $[p]$ at
a vertex of $[P]$ (rather than in its interior) simplifies
calculations, we always try to start with such a choice of $p$
that the system of the level zero roots would be of the maximal
rank, namely $\dim \LL-1$.

If Vinberg's sequence, $v_1,\dots,v_m$, is found up to level $r$,
then the vectors $v\in\Phi$ of higher levels should satisfy the
conditions: $pv<0$ and $vv_i\le0$ for all $v_i$, $1\le i\le m$. If
vectors $v$ respecting these conditions do exist, then the next
segment of Vinberg's sequence is constituted by all such vectors
of the minimal level. Note that the order of Vinberg's roots
within the same level is not well-defined (and is inessential).

This process terminates and gives the basis $\Phi^b$, if the
latter is finite, otherwise the process enumerates vectors of
$\Phi^b$ in an infinite sequence. If we found Vinberg's vectors
$v_1,\dots,v_m$ up to some level $r$, then we can use one of
Vinberg's sufficient criteria below for detecting the termination
of the process.

\subheading{5.4. Vinberg's termination criteria} The Gram matrix
$G_J$ and the Coxeter graph $\Gamma_J$ are called {\it elliptic}
(of rank $r$) if $G_J$ is positive definite (of rank $r$). As is
observed in \cite{Vinberg1}, the elliptic subgraphs of $\G$ of
rank $n-k$ are in one-to-one correspondence with the
$k$-dimensional faces of $[P]$. Namely, an elliptic subgraph
$\G_J$ corresponds to the face supported by the projectivization
of the linear space $H_J=\cap_{_{v\in J}}H_v$.

The connected components of an elliptic graph $G_J$ must belong to
the list of the classical elliptic graphs of the root systems. In
our case (since $m_{vw}=2$ do not appear), an elliptic graph
cannot be other than $A_n$, $D_n$, $E_6,E_7,E_8$,  and $G_2$.

A connected subgraph $\G_J$ and its Gram matrix $G_J$ are called
{\it parabolic} if $G_J$ is positive semi-definite matrix of rank
$|J|-1$. In our case, a parabolic connected subgraph should be one
of the graphs $\widetilde A_n$ (recall that $\tilde A_1$ is just a
thick edge), $\widetilde D_n$, $\widetilde E_6,\widetilde
E_7,\widetilde E_8$, and $\widetilde G_2$, where the subscript
always equals the {\it rank of parabolic graph}, $|J|-1$. A
disconnected subgraph $\G_J$ and its Gram matrix are called
parabolic if all the connected components of $\G_J$ are parabolic.
The rank of such $\G_J$ is by definition the sum of the ranks of
its components. As is observed in \cite{Vinberg1}, a subgraph
$\G_J$ is parabolic of maximal possible rank, $n-1$, if and only
if the intersection $H_J$ defines a vertex of $[P]$ at infinity
(on the Absolute).

Matrix $G_J$ (and its Coxeter graph $\G_J$) is called {\it
critical}, if it is not elliptic, but any submatrix $G_{J'}$,
$J'\varsubsetneq J$, is elliptic. Such $G_J$ is parabolic if
degenerate. If a critical matrix $G_J$ is non-degenerate, its
graph $G_J$ is called {\it Lann\'{e}r's diagrams}.
 The list of Lann\'{e}r's diagrams can be found, for example, in
\cite{Vinberg1}, \cite{Vinberg2}. Note that the only Lann\'{e}r's
diagram possible under the assumptions of this section is a dotted
edge (the other Lann\'{e}r's diagrams all contain a pair of roots
which have the ratio of length different from $1$ and $3$).

\tm{5.4.1. Finite volume criterion \cite{Vinberg1}} Vinberg's
sequence terminates at $J=\{v_1$, $\dots,v_m\}$ if the polyhedron
$P_J$ bounded in $\Lambda_\LL$ by the hyperplanes dual to $v\in J$
has a finite hyperbolic volume.\qed\endtm

To determine finiteness of the volume, Vinberg gives several
criteria. One of them, \cite{Vinberg2} Proposition 4.2(1), can be
formulated (in the form of \cite{Dolgachev} Proposition 2.4) as
follows.

 \tm{5.4.2. Criterion of finiteness of the volume}
The polyhedron $P_J$ has a finite volume if and only if the
following two conditions are satisfied: \roster\item $\G_J$
contains an elliptic subdiagram of rank $n-1$ where $n=\dim\LL-1$,
\item any elliptic subdiagram of rank $n-1$ of $\G_J$ can be
extended to an elliptic subdiagram of rank $n$, or to a parabolic
subdiagram of rank $n-1$; and there exist precisely two such
extensions.\qed
\endroster
\endtm

\rk{Remark} The second condition in theorem 5.4.2 means just that
any edge is adjacent to two vertices: finite, or at infinity.
\endrk

There is another (more simple, but only sufficient) criterion
which can be used if $\G_J$ does not contain Lann\'{e}r's schemes
(that is, dotted edges in our setting).

\tm{5.4.3. Sufficient criterion of finiteness of the volume
\cite{Vinberg1}} The volume of $P_J$ is finite if the following
conditions are satisfied: \roster \item $J$ has rank
$\dim\LL=n+1$;
 \item the Coxeter graph,
$\G_J$, does not contain Lann\'{e}r's diagrams as subgraphs; \item
every connected parabolic subgraph in $\G_J$ is a connected
component of some parabolic subgraph of rank $n-1$ in $\G$.
\endroster
\endtm

\heading \S6. Chirality of $M$-cubics \endheading

\subheading{6.1. Preliminaries and the main statement} A
particular, characteristic, feature of $M$-cubics is that the
lattice $\M$ splits into a direct sum of the eigen-lattices $\M_+$
and $\M_-$. Thus, the eigen-lattices are unimodular in the case of
$M$-cubics, and only in this case. As it follows from the
classification in \cite{FK} (or can be easily deduced directly
from Theorem 4.1.2, Lemma 3.1.1, and the classification of
unimodular lattices), there exists precisely three coarse
deformation classes (equivalently, three homological types) of
$M$-cubics. The corresponding three lattices $\M_+$ are
$U+3I=-I+4I$, $U+3I+E_8=-I+12I$, and $U+3I+2E_8=-I+20I$. The
polarization class $h\in\M_+$ is characteristic, of square $3$,
and can be identified with $(1,1,1)\in 3I$. So, the primitive
lattices $\M_+^0$ are even and isomorphic to $U+A_2$, $U+A_2+E_8$
and $U+A_2+2E_8$ respectively. The corresponding lattices $\M_-$
are also even and isomorphic to $U+2E_8$, $U+E_8$ and $U$
respectively.

\tm{6.1.1. Theorem} Non-singular real cubic fourfold of types
$\M^0_+(c)=U+A_2$ and $\M^0_+(c)=U+A_2+E_8$ are chiral; in
particular, the cubic fourfolds of each of these two types form
two deformation classes. Non-singular real cubic fourfold of type
$\M^0_+(c)=U+A_2+2E_8$ are achiral; these cubic fourfold form one
deformation class.
\endtm

The rest of this Section is devoted to a case by case proof of
this theorem.

We fix a basis $u_1,u_2$ in $U$ and a basis $a_1,a_2$ in $A_2$, so
that $u_i^2=0$ ($i=1,2$), $u_1u_2=1$, $a_i^2=2$ ($i=1,2$), and
$a_1a_2=-1$. The basis $e_1,\dots,e_8$ in $E_8$ is chosen as is
shown on the Coxeter graph of $E_8$, see Figure 1 (we use the
usual convention: $e_i^2=2$ for $i=1,\dots, 8$ and $e_i\circ e_j=
-\delta_{ij}$). This figure presents also the dual vectors
$e_i^*$, $i=1,\dots,8$, which are also elements of $E_8$, because
the lattice $E_8$ is unimodular; for example,
$e_8^*=2e_8+3e_7+4e_6+5e_5+6e_4+4e_3+3e_2+2e_1$. In the case of
$U+A_2+ 2E_8$, the basic vectors of the additional $E_8$-summand
will be denoted by $e_i'$ and their duals by $(e_i')^*$,
$i=1,\dots,8$.

\midinsert \epsfbox{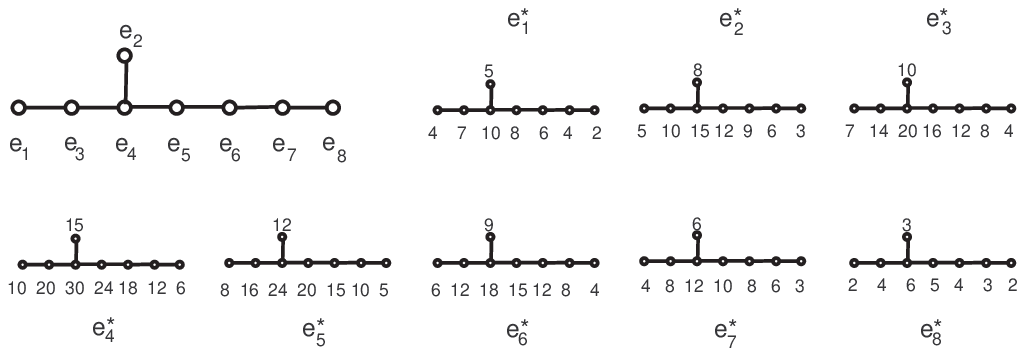} \topcaption{Figure 1. Coxeter graph
$E_8$ and the vectors $e_i^*$}
\endcaption
\endinsert

In all of these three $M$-cases, to apply Vinberg's algorithm (see
Section 5.3) we pick $p=u_1-u_2$. Then we choose as the set of
level-zero vectors the standard bases in each of $E_8$-components
of $\M^0_+(c)$ and complete them by two square-$2$ vectors
$v_1=u_1+u_2$ and $v_2=a_2$, and one square-$6$ vector
$v_3=a_1-a_2$. This choice determines uniquely a cell $P_+$ in
$\L_+(c)$. The vectors of higher levels in Vinberg's sequence must
have components $x_1u_1+x_2u_2+y_1a_1+y_2a_2$ in $U+A_2$
satisfying the following relations:
 $$
 x_2 < x_1,\quad
 x_1+x_2
 \le 0,\quad
 2y_2
 \le y_1,\quad
 y_1
 \le y_2.
 $$
Note that the vector $v_4=-(u_2+a_1+a_2)$ satisfies these relation
and, thus, appears in the list as a vector of level one in each of
the three $M$-cases.

Certainly, the basic vectors of the $E_8$-summands impose also
restrictions on the vectors of higher levels. Namely, their
components in the first (respectively, second) $E_8$-summand
should be linear combinations of $e_1^*,\dots,e_8^*$
(respectively, $(e'_1)^*,\dots,(e'_8)^*$) with non-positive
coefficients.

\subheading{6.2. The case $\M^0_+(c)=U+A_2$} The Coxeter graph of
the vector system $\{v_1,v_2,v_3,v_4\}$ is shown on Figure 2. The
only its parabolic subgraph is $\widetilde{G_2}$ (the subgraph
generated by $v_2$ and $v_3$), and it has rank $2=\dim\L_+-1$. By
Vinberg's finite volume criterion, it implies that Vinberg's
sequence terminates at $\{v_1,v_2,v_3,v_4\}$, and so the
polyhedron $P_+$ is found. Since the Coxeter graph admits no
symmetries, $-\id$ is the only element of $\Aut^-(P_+)$. Thus,
applying Theorem 4.4.1 we conclude that the studied cubic
fourfolds are chiral.

 \midinsert \line{\vtop{\hsize 6cm\vskip-15mm\topcaption{Figure
2}\endcaption \hskip3mm \epsfbox{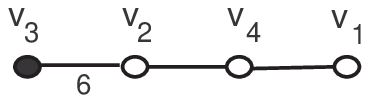} \vskip11.5mm
\centerline{Coxeter's graph for $U+A_2$} }
 \kern2pt\vtop{\hsize 65mm
  $
\matrix\text{Table 1}\\ \\ \boxed{\matrix
&U&A_2\\
\text{------}&\text{------}&\text{------}\\
p&1,-1&0,0\\
\text{level 0}&&\\
v_1&1,1&0,0\\
v_2&0,0&0,1\\
v_3&0,0&1,-1\\
\text{level 1}&&\\
v_4&0,-1&-1,-1\\
\endmatrix}\\ \\
\text{Vinberg's vectors for $\M^0_+(c)=U+A_2$}
\endmatrix
$} }\endinsert

\subheading{6.3. The case $\M^0_+(c)=U+A_2+E_8$} Here, the
level-zero vectors are $e_1,\dots,e_8,$ $ v_1,v_2,$ and $v_3$. The
level-one vectors are $v_4$ and $v_5=-u_2-e_8^*$. This gives the
Coxeter graph shown on Figure 3. This graph has only two parabolic
subgraphs: $\widetilde{G_2}$ and $\widetilde{E_8}$. Vinberg's
finite volume criterion is satisfied because these subgraphs are
disjoint and the sum of their ranks is $2+8=\dim\L_{\M^0_+}-1$.
 The graph has no symmetries and arguing
like in 6.2 we conclude that the studied cubic fourfolds are
chiral.

\midinsert \line{\vtop{\hsize 5cm\vskip-15mm \topcaption{Figure
3}\endcaption \hskip-4mm \epsfbox{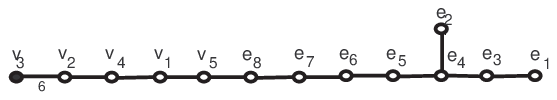} \vskip20.5mm
\centerline{Coxeter's graph for $U+A_2+E_8$} }
 \kern2pt\vtop{\hsize 70mm
$\matrix\text{Table 2
}\\
\\
\boxed{\matrix
&U&A_2&E_8\\
\text{------}&\text{------}&\text{------}&\text{------}\\
p&1,-1&0,0&0\\
\text{level 0}&&&\\
v_1&1,1&0,0&0\\
v_2&0,0&0,1&0\\
v_3&0,0&1,-1&0\\
\text{level 1}&&&\\
v_4&0,-1&-1,-1&0\\
v_5&0,-1&0,0&-e_8^*\\
\endmatrix}\\ \\
\text{Vinberg's vectors for $\M^0_+(c)=U+A_2+E_8$}
\endmatrix
$}}
\endinsert

\subheading{6.4. The case $\M^0_+(c)=U+A_2+2E_8 $} Here, the
level-zero vectors are $e_1,\dots,e_8$, $e_1',\dots,e_8'$,
$v_1,v_2$, and $v_3$. The level-one consists of three $2$-roots
$v_4$, $v_5$, and $v_5'=-u_2-(e_8')^*$. On the next level, $16$,
there is one 2-root $v_6= 2(u_1-u_2)-(a_1+a_2)-e_1^*-(e_1')^*$.
Then, on the level 36 there is a pair of 2-roots:
 $$\aligned
 v_7=3(u_1-u_2)-(2a_1+a_2)-&e_7^*-(e_2')^*,\\
 v_7'=3(u_1-u_2)-(2a_1+a_2)-&e_2^*-(e_7')^*.
\endaligned$$

$$\matrix\text{Table 3. Vinberg's vectors for $\M^0_+(c)=U+A_2+2E_8$}\\
\\
\boxed{\matrix
v_i&U&A_2&E_8&E_8\\
\text{------}&\text{------}&\text{------}&\text{------}&\text{------}\\
p&1,-1&0,0&0&0\\
\text{level 0}&&&&\\
v_1&1,1&0,0&0&0\\
v_2&0,0&0,1&0&0\\
v_3&0,0&1,-1&0&0\\
\text{level 1}&&&&\\
v_4&0,-1&-1,-1&0&0\\
v_5&0,-1&0,0&-e_8^*&0\\
v_5'&0,-1&0,0&0&-(e_8')^*\\
\text{level 16}&&&&\\
v_6&2,-2&-1,-1&-e_1^*&-(e_1')^*\\
\text{level 36}&&&&\\
 v_7&3,-3&-2,-1&-e_7^*&-(e_2')^*\\
 v_7'&3,-3&-2,-1&-e_2^*&-(e_7')^*\\
\text{level 48}&&&&\\
v_8&6,-6&-4,-2&-3e_8^*&-3(e_1')^*\\
v_8'&6,-6&-4,-2&-3e_1^*&-3(e_8')^*\\
\endmatrix}\endmatrix
$$

Our list of Vinberg's vectors given in Table 3 below includes also
a pair of $6$-roots of level $48$,
$$\aligned
v_8=&6(u_1-u_2)-(4a_1+2a_2)-3e_8^*-3(e_1')^*,\\
v_8'=&6(u_1-u_2)-(4a_1+2a_2)-3e_1^*-3(e_8')^*.
\endaligned$$

\midinsert \epsfbox{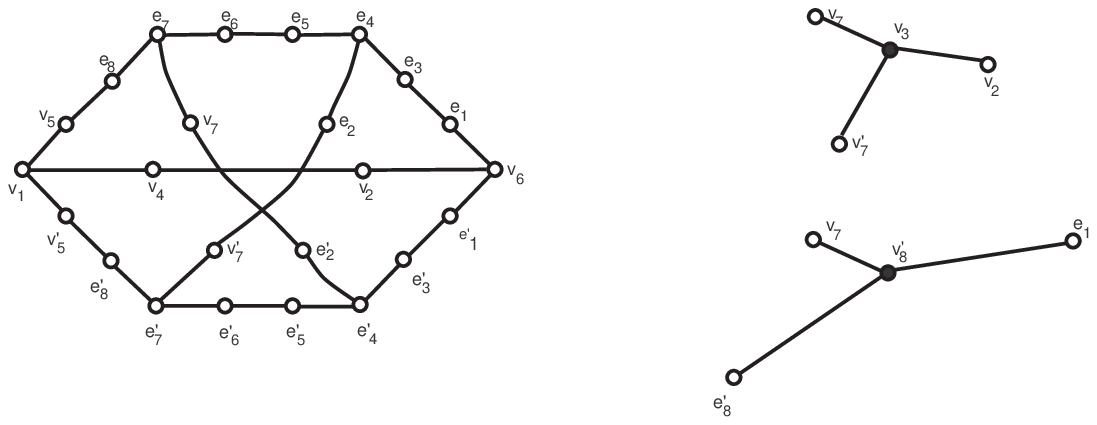} \botcaption{Figure 4. Hexagonal
Coxeter's subgraph for $U+A_2+2E_8$}
\endcaption
\endinsert

The above list contains three 6-roots: $v_3$, $v_8$, and $v_8'$.
If we drop them and consider the Coxeter subgraph formed only by
the 2-roots, we obtain the hexagonal diagram shown on Figure 4.
This diagram has a lot of symmetries. Consider the involution
which fixes the vertices $e_7, e_4',v_6,e_7'$ and permutes the
vertices $v_1, e_4$. Since the set of vectors corresponding to the
vertices of the diagram generate the lattice $\M^0_+(c)$, this
involution is induced by a lattice involution $f:\M^0_+(c)\to
\M^0_+(c)$ (see Proposition 5.2.1). Since in the whole Coxeter
diagram the 6-root $v_3$ is connected with the 2-roots
$v_2,v_7,v'_7$ and the 6-root $v_8'$ is connected with the 2-roots
$e_1,v_7,e_8'$, the automorphism $f$ transforms $v_3$ into $v_8'$.
The $A_2$-components of $v_3$ and $v_8'$ are $(1,-1)$ and
$(-4,-2)$, which are not congruent modulo $3$. This implies that
$f$ is $\Z/3$-reversing. By Proposition 5.2.1, $f$ is
$P_+$-direct, so applying Theorem 4.4.1 we conclude that this
homological type is achiral.

\rk{Remark} This lattice, its fundamental chamber, and the
complete Coxeter graph had appeared already in Vinberg's paper
\cite{V3} on maximally algebraic $K3$-surfaces. Note that our list
contains the full set of 2-roots, and the missing 6-roots can be
obtained from the 6-roots in the list by applying the symmetries
of the hexagonal subgraph. The same construction is given in
\cite{Looijenga}.
\endrk

\heading \S7. Chirality of $(M-1)$-cubics\endheading

\subheading{7.1. Preliminaries and the main statement} The next
after $M$-cubics by their topological complexity are
$(M-1)$-cubics. The deformation components of the latters, as it
follows from \cite{FK}, are adjacent to the deformation components
of $M$-cubics. The lattice, $\M$, of an $(M-1)$-cubic contains the
direct sum of the eigen-lattices $\M_+$ and $\M_-$ as a sublattice
of index $2$, and this condition characterizes $(M-1)$-cubics
among all non-singular real cubic fourfolds. In the other words,
the characteristic feature of $(M-1)$-cubics is that $\M_\pm$ have
discriminant $\Z/2$. Using the general properties of lattices
$\M_\pm$ (namely, that lattice $\M_+$ is odd with a characteristic
element $h\in \M_+$ of square $h^2=3$, that lattice $\M_-$ is
even, and that the both lattices are of index $\sigma_-=1$), one
can deduce that the $(M-1)$-cubics form precisely six homological
types, see \cite{FK}. As usual, these types can be distinguished
by sublattices $M_+$, as well as by sublattices $M_-$. The
corresponding six lattices $\M^0_+$ are $U+A_2+A_1+kE_8$ and
$-A_1+A_2+kE_8$, $k=0,1,2$.

\tm{7.1.1. Theorem} Non-singular real cubic fourfolds of types
$\M^0_+(c)=-A_1+A_2, U+A_2+A_1$,  and $-A_1+A_2+E_8$, are chiral;
in particular, the cubic fourfolds of each of these three types
form two deformation classes. Non-singular real cubic fourfolds of
types $\M^0_+(c)=U+A_2+E_8+A_1, -A_1+A_2+2E_8$, and
$U+A_2+2E_8+A_1$ are achiral; the cubic fourfolds of each of these
three types form one deformation class.
\endtm

\subheading{7.2. The case $\M_+^0(c)=-A_1+A_2$} Here, Vinberg's
sequence starts from vectors $\{v_1,v_2,v_3\}$ given in Table 4.
The Coxeter graph of this sequence of three vectors is shown in
Figure 5. It contains a unique parabolic subgraph
$\widetilde{A_1}$ (a thick edge connecting $v_2$ and $v_3$).
Vinberg's criteria 5.4.3 and 5.4.1 can be applied to conclude
termination, since the rank of $\widetilde A_1$ is $1=\dim
\M_+^0(c)-2$. The Coxeter graph admits no symmetries. Hence,
applying Theorem 4.4.1 we deduce that the studied cubic fourfolds
are chiral.

 \midinsert
 \line{\vtop{\hsize 5cm\vskip-20.75mm\topcaption{Figure
5}\endcaption \hskip-3mm\vskip2mm \epsfbox{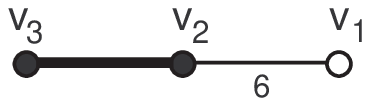}\vskip8mm
 \centerline{Coxeter's
graph for $-A_1+A_2$} }
 \kern2pt\vtop{\hsize 70mm
$\matrix\text{Table 4}\\ \\ \boxed{\matrix
&-A_1&A_2\\
\text{------}&\text{------}&\text{------}\\
p&1&0,0\\
\text{level 0}&&\\
v_1&0&0,1\\
v_2&0&1,-1\\
\text{level 12}&&\\
v_3&3&-4,-2\\
\endmatrix}\\ \\ \text{Vinberg's vectors for $\M^0_+(c)=-A_1+A_2$}
\endmatrix
$}}\endinsert

\subheading{7.3. The case $\M_+^0(c)=U+A_2+A_1$} Here, Vinberg's
sequence starts from four level-zero vectors $\{v_1,v_2,v_3,
v_4\}$  and two level-one vectors $\{v_5,v_6\}$ given in Table 5.
The Coxeter graph of this sequence of six vectors is shown in
Figure 6. It contains precisely two parabolic subgraphs,
$\widetilde G_2$ (vertices $v_3$, $v_2$, $v_5$) and $\widetilde
A_1$ ($v_4$, $v_6$). Vinberg's criterion is satisfied, since the
rank of their union is $2+1=\dim \M_+^0(c)-2$. The Coxeter graph
admits no symmetries. Hence, applying Theorem 4.4.1 we conclude
that the studied cubic fourfolds are chiral.

 \midinsert \line{\vtop{\hsize 5.5cm\vskip-25.75mm\topcaption{Figure
6}\endcaption \hskip-6mm\vskip10mm \epsfbox{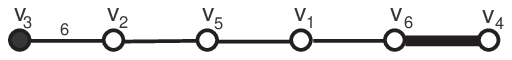} \vskip23mm
\centerline{Coxeter's graph for $U+A_2+A_1$} }
 \kern2pt\vtop{\hsize 75mm\hskip-3mm
$\matrix\text{Table 5}
\\ \\
\boxed{\matrix
&U&A_2&A_1\\
\text{------}&\text{------}&\text{------}&\text{------}\\
p&1,-1&0,0&0\\
\text{level 0}&&&\\
v_1&1,1&0,0&0\\
v_2&0,0&0,1&0\\
v_3&0,0&1,-1&0\\
v_4&0,0&0,0&1\\
\text{level 1}&&&\\
v_5&0,-1&-1,-1&0\\
v_6&0,-1&0,0&-1\\
\endmatrix}\\ \\ \text{Vinberg's vectors for $\M^0_+(c)=U+A_2+A_1$}
\endmatrix
$}}\endinsert

\subheading{7.4. The case $\M_+^0(c)=-A_1+A_2+E_8$} Here, the
level-zero vectors of Vinberg's sequence are $e_1,\dots,e_8$,
$v_1$, and $v_2$. They are followed by two vectors of level 4 and
one vector of level 12, see Table 6. The Coxeter graph, $\G$, of
this sequence of thirteen vectors is shown in Figure 7.

  \midinsert\line{\vtop{\hsize 5cm\vskip-27.25mm\topcaption{Figure
7}\endcaption \hskip-3mm \vskip5mm\epsfbox{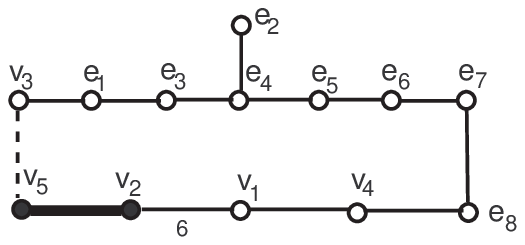} \vskip10mm
\centerline{Coxeter's graph for $-A_1+A_2+E_8$} }
 \kern2pt\vtop{\hsize 75mm
$\matrix\text{Table 6}\\
\\
\boxed{\matrix
&-A_1&A_2&E_8\\
\text{------}&\text{------}&\text{------}&\text{------}\\
p&1&0,0&0\\
\text{level 0}&&\\
v_1&0&0,1&0\\
v_2&0&1,-1&0\\
\text{level 4}&&&\\
v_3&1&0,0&-e_1^*\\
v_4&1&-1,-1&-e_8^*\\
\text{level 12}&&&\\
v_5&3&-4,-2&0\\
\endmatrix}\\ \\ \text{Vinberg's vectors for $\M^0_+(c)=-A_1+A_2+E_8$}
\endmatrix
$}}\endinsert

\tm{Lemma 7.4.1} Vinberg's criterion 5.4.2 is satisfied for the
Coxeter graph $\G$ on Figure 7.
\endtm

\demo{Proof} For $S=\{a_1,\dots,a_n\}\subset\Phi^b$, let
$F_S=F_{a_1,\dots,a_n}\subset P$ denote the face of the cell $P$
supported in the intersection of the walls $[H_v]$, where
$v\in\Phi^b\smallsetminus S$. Note that $P$ has two vertices at
infinity, $F_{v_5,e_8}$ and $F_{v_1, v_3}$ (because the sets
$\Phi^b\smallsetminus S$ span parabolic subgraphs of maximal
possible rank $\dim(\M_+^0)-2=9$). The other vertices of $P$ are
$F_{a,b,c}$ such that $\Phi^b\smallsetminus\{a,b,c\}$ spans an
elliptic subgraph. This subgraph cannot contain the dotted edge
connecting $v_3$ with $v_5$, so the set $S=\{a,b,c\}$ should
contain either $v_3$ or $v_5$ (or the both). This set should
contain also at least one vertex-root from each of the parabolic
subgraphs $\widetilde E_7$, $\widetilde E_8$, $\widetilde G_2$,
$\widetilde A_1$ of $\G$. If the both $v_3$ and $v_5$ are included
in $S=\{a,b,c\}$, then $F_S$ is a vertex of $P$ only for
$S=\{v_3,v_4,v_5\}$. If $a=v_5$ and $v_3\notin S$, then $b$ and
$c$ should be chosen from the two disjoint parabolic subgraphs
$\widetilde{G_2}$ and $\widetilde{E_7}$, which gives $21$ other
vertices $F_{v_5,b,c}$, where $b\in\{v_2,v_1,v_4\}$ and
$c\in\{e_1,\dots, e_7\}$. Similarly, if $a=v_3$ and $v_5\notin S$,
then $b$, $c$ should be chosen from the two disjoint parabolic
subgraphs $\widetilde{A_1}$ and $\widetilde{E_8}$, so $b=v_2$ and
$c\in\{e_1,\dots,e_8,v_4\}$, which gives $9$ new vertices
$F_{a,b,c}$. Totally, $\G$ contains $31$ finite vertices and two
vertices at infinity.

The edges of $P$ can be expressed as $F_S$, $S=\{a,b,c,d\}$, where
$\Phi^b\smallsetminus S$ spans an elliptic subgraph. Thus, as
above, $S$ should contain at least one of $v_3$ and $v_5$.
 In the edges $F_{v_3,v_5,v_4,d}$, one of the endpoints is
 $F_{v_3,v_5,v_4}$. In the cases $d\in\{e_1,\dots,e_7\}$, the
 other endpoint is $F_{v_5,v_4,d}$. In the cases $d=v_1$, $d=v_2$,
 and
 $d=e_8$, the other endpoint is respectively $F_{v_1,v_3}$,
 $F_{v_3,v_4,v_2}$,
 and
 $F_{v_5,e_8}$.
The edges $F_{v_3,v_5,e_8,d}$ must have $d\in\{v_1,v_2\}$ and are
incident to $F_{v_5,e_8}$. Another endpoint is $F_{v_3,v_1}$ for
$d=v_1$, and $F_{v_3,e_8,v_2}$ for $d=v_2$. Each of the edges
$F_{v_3,v_5,c,d}$, $c\in\{v_1,v_2\}$, $d\in\{e_1,\dots,e_7\}$ has
$F_{v_5,c,d}$ as one of the endpoints. The other endpoint is
$F_{v_3,v_2,d}$ if $c=v_2$ and $F_{v_3,v_1}$ if $c=v_1$.

The other edges $F_{a,b,c,d}$ have $a\in\{v_3,v_5\}$ and
$b,c,d\notin\{v_3,v_5\}$. If $a=v_3$, then another vertex should
be chosen from the subgraph $\widetilde A_1$, and we may assume
that $b=v_2$ (since the case $b=v_5$ was already considered). This
gives edges $F_{v_3,v_2,c,d}$ with $c\in\{e_1,\dots,e_8,v_4\}$ and
$d\in\{e_1,\dots,e_8,v_4,v_1\}$. If $d\ne v_1$, then the endpoints
are $F_{v_3,v_2,c}$ and $F_{v_3,v_2,d}$. The endpoints of
$F_{v_3,v_2,c,v_1}$ are $F_{v_3,v_2,c}$ and $F_{v_3,v_1}$.
 Finally, if $a=v_5$, then one of $b,c,d$ should
be chosen from $\widetilde G_2$, say, $b\in\{v_2,v_1,v_4\}$ and
another from $\widetilde E_7$, say, $c\in\{e_1,\dots,e_7\}$.
 Then $F_{v_5,b,c,d}$ has one endpoint $F_{v_5,b,c}$. Another
its endpoint is $F_{v_5,c,d}$ if $d\in\{v_2,v_1,v_4\}$, and
$F_{v_5,b,d}$ if $d\in\{e_1,\dots,e_7\}$. In the remaining case
$d=e_8$, the second endpoint is $F_{v_5,e_8}$. \qed\enddemo

The Coxeter graph admits no symmetries. Hence, applying Theorem
4.4.1 we conclude that the studied cubic fourfolds are chiral.

\subheading{7.5. The case $\M_+^0(c)=U+A_2 +A_1 +E_8$} Here, the
level-zero Vinberg's vectors are $e_1,\dots,e_8$ plus $v_1,\dots,
v_4$ listed in Table 7. Then follow three vectors $v_5,v_6,v_7$ of
level 1 and the vector $v_8$ of level $48$ (see the same Table).
\vskip2mm

  \midinsert\line{\vtop{\hsize 5cm\vskip-24.5mm\topcaption{Figure
8}\endcaption \hskip-3mm\vskip10mm \epsfbox{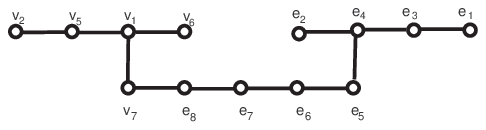} \vskip13mm
\centerline{A symmetric fragment}\centerline{of Coxeter's graph}
\centerline {for $U+A_2+A_1+E_8$} }
 \kern2pt\vtop{\hsize 75mm\hskip-7mm
$\matrix\text{Table 7}\\
\\
\boxed{ \matrix
&U&A_2&A_1&E_8\\
\text{------}&\text{------}&\text{------}&\text{------}&\text{------}\\
p&1,-1&0,0&0&0\\
\text{level 0}&&&&\\
v_1&1,1&0,0&0&0\\
v_2&0,0&0,1&0&0\\
v_3&0,0&1,-1&0&0\\
v_4&0,0&0,0&1&0\\
\text{level 1}&&&&\\
v_5&0,-1&-1,-1&0&0\\
v_6&0,-1&0,0&-1&0\\
v_7&0,-1&0,0&0&-e_8^*\\
\text{level 48}&&&&\\
v_8&6,-6&-4,-2&-3&-3(e_1)^*.
\endmatrix}\\ \\ \text{Vinberg's vectors for $\M^0_+(c)=U+A_2+A_1+E_8$}
\endmatrix
$}}\endinsert

Consider the Coxeter subgraph formed by Vinberg's vectors
$e_1,\dots, e_8$, $v_1, v_2, v_5$, $v_6, v_7$. This subgraph is
shown in Figure 8. It has an evident nontrivial involution (which
fixes the vertex $v_7$ and permutes the vertices $v_2, e_1$).
Since the vectors $e_1, \dots, e_8$, $v_1, v_2, v_5, v_6, v_7$
generating this subgraph span the lattice $\M^0_+(c)$, this
involution is induced by a $P_+$-direct lattice involution
$f:\M^0_+(c)\to \M^0_+(c)$ (see Proposition 5.2.1). In particular,
$f$ transforms Vinberg's vector $v_3=-v_5-2v_2+v_7+e_8^*$, into
another Vinberg's vector
$v'_3=-e_3-2e_1+e_5+(2e_6+3e_7+4e_8+5v_7+6v_1+4v_5+3v_6+2v_2)$.
The $A_2$-component of $v'_3$ is $4(-1,-1)+2(0,1)=(-4,-2)$, while
the $A_2$-component of $v_3$ is $(1,-1)$. Hence, $f$ is
$\Z/3$-reversing and applying Theorem 4.4.1 we conclude that the
type considered is achiral.

\comment
 \midinsert
\hskip5mm \epsfbox{m-1b2.eps} \botcaption{Figure 8. A symmetric
fragment of Coxeter's graph for $U+A_2+A_1+E_8$}
\endcaption
\endinsert

$$\matrix\text{Table 7. Vinberg's vectors for $\M^0_+(c)=U+A_2+A_1+E_8$}\\
\\
\boxed{ \matrix
&U&A_2&A_1&E_8\\
\text{------}&\text{------}&\text{------}&\text{------}&\text{------}\\
p&1,-1&0,0&0&0\\
\text{level 0}&&&&\\
v_1&1,1&0,0&0&0\\
v_2&0,0&0,1&0&0\\
v_3&0,0&1,-1&0&0\\
v_4&0,0&0,0&1&0\\
\text{level 1}&&&&\\
v_5&0,-1&-1,-1&0&0\\
v_6&0,-1&0,0&-1&0\\
v_7&0,-1&0,0&0&-e_8^*\\
\text{level 48}&&&&\\
v_8&6,-6&-4,-2&-3&-3(e_1)^*.
\endmatrix
}\endmatrix$$
\endcomment

\subheading{7.6. The case $\M_+^0(c)=-A_1+A_2+2E_8$} Let us start
with a bit more general setting. Namely, assume that $\LL$ is a
lattice like in \S5 (for example, some of the lattices
$\M_+^0(c)$), $P\subset\L(\LL)$ is a cell, and $f\in\Aut^+(P)$
 is an automorphism of $\LL$ induced by some
symmetry of the Coxeter graph, $\G$, of $\LL$. Suppose that a
$2$-root $v$ is a vertex of $\G$ preserved by this symmetry, that
is, $f(v)=v$. Then the sublattice $\LL_v=\{x\in\LL\,|\,xv=0\}$ is
$f$-invariant and we may consider an induced automorphism
$f_v\in\Aut(\LL_v)$.

\tm{7.6.1. Lemma} If $f$ is $\Z/3$-reversing, then $f_v$ is also
$\Z/3$-reversing, and $P_v$-direct for some cell $P_v$ of
$\LL_v$.\endtm

\demo{Proof} Since $\discr_3(\LL_v)=\discr_3(\LL)=\Z/3$, the
automorphisms $f$ and $f_v$ are the both $\Z/3$-direct or
$\Z/3$-reversing. Furthermore, $f$ preserves the facet
$P\cap[H_v]$ of $P$, since it preserves both $P$ and $v$. Due to
$\discr_3(\LL_v)=\discr_3(\LL)$, each wall in $\L(\LL_v)$ is an
intersection of $[H_v]$ with a wall $\L(\LL)$, and thus, the facet
$P\cap[H_v]$ is a part of some cell, $P_v$, of $\L(\LL_v)$. Such
$P_v$  has to be also invariant.
 \qed\enddemo

\tm{7.6.2. Corollary} Lattice $-A_1+A_2+2E_8$ is achiral.
\endtm

\demo{Proof} Let $\LL=U+A_2+2E_8$, then for $v=v_1$ (notation of
 6.4) we have $\LL_v=-A_1+A_2+2E_8$.
An involution of the hexagonal diagram (Figure 4) in 6.4 is
conjugate to some involution, $f\in\Aut^+(P)$, preserving $v_1$.
Since $f$ is $\Z/3$-reversing, we can apply Lemma 7.6.1.
\qed\enddemo

Applying Theorem 4.4.1 we can now conclude that the fourfolds with
$\M_+^0(c)=-A_1+A_2+2E_8$ are achiral.

\subheading{7.7. The case $\M^0_+(c)=U+A_2+2E_8+A_1$} Let $\LL$
and $\LL_v$ be like in 7.6. Our aim now is to obtain a criterion
which is in some sense ``converse'' to the one in Lemma 7.6.1.
Recall that lattice $\LL$ either splits into a direct sum of
$\LL_v$ with a sublattice $A_1=\Z v$, or contains this direct sum
as an index $2$ sublattice. We will show that, in the former case,
achirality of $\LL_v$ implies achirality of $\LL$.

\tm{7.7.1. Lemma} Assume that $\LL=\LL_v+A_1$, where $A_1=\Z v$,
and $\LL_v$ is achiral. Then $\LL$ is also achiral. In fact, any
$\Z/3$-reversing automorphism $f_v\in\Aut^+(P_v)$ for some cell
$P_v\subset\L(\LL_v)$ can be extended to a $\Z/3$-reversing
automorphism $f\in\Aut^+(P)$ for some cell $P\subset\L(\LL)$.
\endtm

\demo{Proof} Letting $f(v)=v$ we obtain an extension of $f_v$ to
$\LL$ which is obviously $\Z/3$-reversing if $f_v$ is.

Like in Lemma 7.6.1, by the same evident reasons, $P_v$ contains
the facet $P\cap[H_v]$ of some cell $P$ in $\L(\LL)$. But now, the
relation is stronger: $P\cap[H_v]=P_v$. In fact, the walls of $P$
different from $[H_v]$ are either orthogonal to $[H_v]$ or do not
intersect it. To see it, consider any wall $[H_w]$, $w\in V_2\cup
V_6$. Splitting $\LL=\LL_v+\Z v$ gives a decomposition $w=w_v+kv$,
where $w_v\in\LL_v$, $k\in\Z$. If $k=0$, then $[H_w]$ is
orthogonal to $[H_v]$, whereas $w_v=0$ implies $w=v$. Otherwise we
observe that $w_v^2=w^2-k^2v^2\le0$, because $v^2=2$, and $w^2$ is
either $2$, or $6$, but in the latter case $k$ is divisible by
$3$. Thus, vectors perpendicular to $w_v$ cannot have negative
square, which contradicts to $P\cap[H_v]\ne\emptyset$.

The relation $P\cap[H_v]=P_v$ implies that the isometry
$f=f_v\oplus\id:\LL\to \LL$  is $P$-direct.
 \qed\enddemo

\tm{7.7.2. Corollary} Lattice $U+A_2+2E_8+A_1$ is achiral.
\endtm

\demo{Proof} According to 6.4, the lattice $\LL_v=U+A_2+2E_8$ is
achiral. It remains to apply Lemma 7.7.1. \qed\enddemo

Applying Theorem 4.4.1 we can now conclude that the cubic
fourfolds with $\M^0_+(c)=U+A_2+2E_8+A_1$ are achiral.

\heading \S8. Concluding remarks \endheading

\subheading{8.1. Further results}
 The cases of $M$-varieties and $(M-1)$-varieties are
usually the most interesting and difficult ones, which explains
our special interest to them in the context of the chirality
problem of the cubic fourfolds. But our methods are applicable as
well to the other cases. Our observations concerning the problem
of chirality can be summarized as follows.

 Let $\rho$ denote
the rank of the lattice $\M_+^0$, $r=22-\rho$ denote the rank of
$\M_-$, and $d$ the discriminant rank,
$\roman{rk}(\discr_2(\M_+^0))=\roman{rk}(\discr(\M_-))$.
 In all the cases studied,
 if $\rho+d\ge 14$
 then  $\M_+^0$ is
 achiral. In addition, the list of achiral lattices contains
 $\M_+^0=U(2)+A_2+D_4$ and
$\M_+^0=-A_1+\la6\ra+kA_1$ with $k=2,3,$ and $4$. The other
lattices that we have analyzed are chiral. (In a few cases
remaining for analysis, the discriminant form is even and
$\rho+d\ge 14$. We expect that the corresponding lattices are
achiral.)

The lattices $M_+^0(c)$ of cubic fourfolds can be naturally
divided into the {\it principal series}, which contains the most
of lattices and is presented in Table 8, and several additional
lattices presented in Table 9 (see \cite{FK} for more details).

\midinsert \line{\vtop{\hsize 7cm
$\matrix\text{Table 8}\\ \text{The principal series of lattices $M_+^0$}\\
\boxed{\matrix
 -A_1&+\la6\ra&+t A_1,\quad 0\le t\le9\\
-A_1&+A_2&+t A_1,\quad 0\le t\le9\\
 U&+A_2&+t A_1,\quad 0\le t\le9\\
 U&+A_2+D_4&+t A_1,\quad 0\le t\le6\\
 -A_1&+\la6\ra+E_8&+t A_1,\quad 0\le t\le5\\
 -A_1&+A_2+E_8&+t A_1,\quad 0\le t\le5\\
 U&+A_2+E_8&+t A_1,\quad 0\le t\le5\\
 U&+A_2+D_4+E_8&+t A_1,\quad 0\le t\le2\\
 -A_1&+\la6\ra+2E_8&+t A_1,\quad 0\le t\le1\\
 -A_1&+A_2+2E_8&+t A_1,\quad 0\le t\le1\\
 U&+A_2+2E_8&+t A_1,\quad 0\le t\le1\\
\endmatrix
}\endmatrix$}\hskip4mm
 \kern2pt\vtop{\hsize 6cm
$\matrix\text{Table 9}\\ \text{Additional lattices $M_+^0$}\\
\boxed{\matrix
U(2)+E_6(2)\\
U(2)+A_2\\
U+E_6(2)\\
 U(2)+A_2+D_4\\
 U(2)+A_2+2D_4\\
 U+A_2+2D_4\\
 U(2)+A_2+E_8\\
 U(2)+A_2+D_4+E_8\\
 U(2)+A_2+2E_8\\
 U+A_2+E_8(2)\\
 U(2)+A_2+E_8(2)\\
\endmatrix}\endmatrix$}}\endinsert

The Table 10 describes chirality of the principal series of cubic
fourfolds in terms of the ranks $r$ and $d$.
\midinsert  \topcaption{Table 10. Chirality of cubic fourfolds:
the principal series}\endcaption \centerline{ $ \smallmatrix
\bold{d}\\
{11}& &&&&&&&&&&a&&&&&&&&&&\\
{10}& &&&&&&&&&a&&a&&&&&&&&\\
9&  &&&&&&&&a&&a&&a&&&&&&\\
8&  &&&&&&&a&&a&&a&&a&&&&&\\
7&  &&&&&&a&&a&&a&&a&&a&&&&\\
6&  &&&&&a&&a&&a&&a&&a&&a&&&\\
5&  &&&&a&&a&&a&&a&&a&&c&&a&&\\
4&  &&&a&&a&&a&&a&&a&&c&&c&&a&&\\
3&  &&a&&a&&a&&a&&a&&c&&c&&c&&c&\\
2&  &a&&a&&a&&a&&a&&c&&c&&c&&c&&c\\
1&  a&&a&&&&&&a&&c&&&&&&c&&c\\
0&  &a&&&&&&&&c&&&&&&&&c&\\
&\,1\,&\,2\,&\,3\,&\,4\,&\,5\,&\,6\,&\,7\,&\,8\,&\,9\,&{10}&{11}&{12}&{13}&{14}&{15}&{16}&{17}&{18}&{19}&{20}&\bold{r}
\endsmallmatrix$}\vskip2mm
\centerline{\eightrm The symbol ``c'' stands for the chiral
deformation classes, and symbol ``a'' for the achiral
ones}\endinsert

\comment
$$\matrix
\bold{d}\\
11&&&&&&&&&&a&&&&&&&&&&&\\
10&&&&&&&&&a&&a(?)&&&&&&&&&&\\
9&&&&&&&&a&&a&&a&&&&&&&\\
8&&&&&&&a(a)&&a&&a(?)&&a&&&&&&\\
7&&&&&&a&&a&&a&&a&&a&&&&&\\
6&&&&&a&&a(a)&&a&&a(?)&&a&&a&&&&\\
5&&&&a&&?&&a&&a&&a&&a&&a&&&\\
4&&&a&&c&&c(a)&&a&&a(a)&&a&&a(?)&&a&&&\\
3&&c&&c&&c&&c&&a&&a&&a&&a&&a&&\\
2&c&&c(c)&&c&&c&&c&&a(a)&&a&&a&&a&&a(?)&\\
1&&c&&c&&&&&&c&&a&&&&&&a&&a\\
0&&&c&&&&&&&&c&&&&&&&&a&\\
\\
&2&3&4&5&6&7&8&9&10&11&12&13&14&15&16&17&18&19&20&21&22&\bold{r}\\
\endmatrix$$
\endcomment

\subheading{8.2. Chirality of singular cubic fourfolds} Chirality
of cubic fourfolds having nodal singularities is an interesting
related problem. It is trivial to observe that any perturbation of
an achiral nodal cubic provides an achiral non-singular cubic. The
non-trivial part of the problem is the converse: if perturbations
give only achiral cubics, can we conclude that a nodal cubic is
achiral itself ? To solve it, one can use the same approach as in
the non-singular case, just taking into account the vanishing
cycles. On the other hand, the central projective correspondence
discussed in \cite{FK} relates chirality of nodal cubic fourfolds
to a certain question about $6$-polarized $K3$-surfaces. This
relation can be used in the both directions.

A somewhat different kind of observation is chirality of the {\it
discriminant cubic}, $\det\left(\matrix
x_0&x_1&x_2\\
x_1&x_3&x_4\\
x_2&x_4&x_5 \endmatrix\right)=0$, which parameterizes the space of
singular conics in the plane. The key observation is chirality of
the singular locus of the discriminant cubic (this locus is the
image of the Veronese map).

\subheading{8.3. Explicit equations} It would be interesting to
find explicit (natural) equations for representatives of each of
the deformation classes. It can be helpful not only for proving
achirality statements, but also for better understanding of the
topology of the cubic hypersurfaces. As an example, let us
consider the equations of the following type:
$$
(\sum^6_1 x_\alpha)^3-\sum^6_1 c_\alpha x_\alpha^3=0;
$$
these equations were proposed in the late 70th by D.~Fucks
(private communication to the second author) for searching the
precise range of the values of the Euler characteristic of real
cubic hypersurfaces in each given even dimension (a problem
remaining, up to our knowledge still open in its whole
generality). Similar equations were used earlier by F.~Klein
\cite{Kl}, and his student C.~Rodenberg \cite{R}, to find and to
study explicit representatives for each of the five classes of
real nonsingular cubic surfaces. In fact, it is by means of these
equations that Klein proved in \cite{Kl2} the achirality of all
real nonsingular cubic surfaces (cf., the remark at the end of
this subsection).

One can easily check that for $c_\alpha$ having all the same value
$c$, the topology of the hypersurface is changing at $c=0, 4, 16$,
and $36$. For $c<0$ and $c>36$ the real part of the hypersurface
is diffeomorphic to the real four-dimensional projective space,
$\Rp4$. When $c=36$, there appears a solitary double point, so
that for $16<c<36$ we observe $S^4\sqcup{\Rp4}$. When $c=16$, our
hypersurfaces acquire six double points of Morse index $(1,4)$
with respect to growing $c$ (the first, respectively second,
component of the index is the number of positive, respectively,
negative squares), and therefore, for $4<c<16$ the real part of
the hypersurface is diffeomorphic to the real four-dimensional
projective space with five $S^1\times S^3$-handles, that is
$\Rp4\#5(S^1\times S^3)$. Finally, when $c=4$, one finds that
there are fifteen double points of Morse index $(2,3)$, which
implies that the Euler characteristic of our hypersurfaces becomes
equal to $21$. According to the classification of cubics (see
\cite{FK}), there is only one coarse deformation class with this
value of Euler characteristic (in fact, it is the class studied
above in Section 7.7), and for the cubics of this class the real
part has the homological type of $\Rp4\#10(S^2\times S^2)$. (One
can also give a direct proof based on the Lefschetz trace formula
and the Smith theory, which allow to reconstruct the Betti numbers
from the action of the complex conjugation in homology.)

Since for $c_\alpha$ having all the same value the equation is
invariant under transposition of the variables, all these
hypersurfaces represent achiral classes. In the same manner, one
can show that the whole left-hand slanted border of the diagram
shown in Table 10 consists exclusively of achiral classes.

\subheading{8.4. Chirality in lower dimensions. Quartic surfaces}
Speaking on the real non-singular hypersurfaces $X_\R$ of
dimension $n$ and degree $d$, it is easy to see their achirality
in the trivial cases $n=0$ (for any $d$), and $d\le2$ (for any
$n$). As was pointed out in Introduction, $X_\R$ is also achiral
if $n$ is odd. Achirality of cubic surfaces was observed by
F.~Klein, as we mentioned in 8.3. The next case of quartic
surfaces was analyzed in \cite{Kh1},\cite{Kh2} using a technique
similar to the one in this paper. It turned out that a
non-contractible (in $\Rp3$) quartic $X_\R$ is chiral if and only
if it has at least 4 spherical components, and a contractible
quartic is chiral if and only if it has at least 3 spherical
components and, in addition, a component with at least 3 handles
(see Table 11, where $r$ is the rank of the $+1$-eigen-lattice
$L_+=\{x\in H_2(X)\: \conj_*x=x\} $, $d$ is the discriminant rank
of $L_+$, and symbols $a$, or $c$ stand as in Table 10 for
achiral, or respectively, chiral deformation classes).

\midinsert  \topcaption{Table 11. Chirality of quartic
surfaces}\endcaption \centerline{ $ \smallmatrix
\bold{d}\\
{10}& &&&&&&&&&a&&&&&&&&&\\
9&  &&&&&&&&a&&a&&&&&&&\\
8&  &&&&&&&a&&a&&a&&&&&&\\
7&  &&&&&&a&&a&&a&&a&&&&&\\
6&  &&&&&a&&a&&a&&a&&c&&&&\\
5&  &&&&a&&a&&a&&a&&c&&c&&&\\
4&  &&&a&&a&&a&&a&&c&&c&&c&&\\
3&  &&a&&a&&a&&a&&c&&c&&c&&c&\\
2&  &a&&a&&a&&a&&c&&c&&c&&c&&c&\\
1&  a&&a&&&&&&c&&c&&&&&&c&&c\\
0&  &a&&&&&&&&c&&&&&&&&c&\\
&\,1\,&\,2\,&\,3\,&\,4\,&\,5\,&\,6\,&\,7\,&\,8\,&\,9\,&{10}&{11}&{12}&{13}&{14}&{15}&{16}&{17}&{18}&{19}&&\bold{r}
\endsmallmatrix$
} \vskip2mm\centerline{\eightrm Non-contractible case} \vskip4mm
\centerline{ $ \smallmatrix
\bold{d}\\
{11}& &&&&&&&&&&a&&&&&&&&&&\\
{10}& &&&&&&&&&a&&a&&&&&&&&\\
9&  &&&&&&&&a&&a&&a&&&&&&\\
8&  &&&&&&&a&&a&&a&&a&&&&&\\
7&  &&&&&&a&&a&&a&&a&&a&&&&\\
6&  &&&&&a&&a&&a&&a&&a&&a&&&\\
5&  &&&&a&&a&&a&&c&&a&&a&&a&&\\
4&  &&&a&&a&&a&&c&&c&&a&&a&&a&&\\
3&  &&a&&a&&a&&c&&c&&c&&a&&a&&a&\\
2&  &a&&a&&a&&c&&c&&c&&c&&a&&a&&a\\
1&  a&&a&&&&&&c&&c&&&&&&a&&a\\
0&  &a&&&&&&&&c&&&&&&&&a&\\
&\,1\,&\,2\,&\,3\,&\,4\,&\,5\,&\,6\,&\,7\,&\,8\,&\,9\,&{10}&{11}&{12}&{13}&{14}&{15}&{16}&{17}&{18}&{19}&{20}&\bold{r}
\endsmallmatrix$
} \vskip2mm \centerline {\eightrm Contractible case}
\endinsert

\subheading{8.5. Reversibility} In connection with chirality, it
may be worth mentioning a different but somehow related notion of
{\it reversibility}, which plays a non-trivial role for instance
for odd-dimensional hypersurfaces. Namely, to each deformation
class of real non-singular hypersurfaces $X\subset\Rp{n+1}$ of
degree $d$, that is a connected component $\Cal C$ of $\Cal
C_{n,d}=P_{n,d}(\R) \setminus\D_{n,g}(\R)$, we can associate its
pull back $\til {\Cal C}$ into the sphere $\til P_{n,d}(\R)$ which
covers $P_{n,d}(\R)$. This pull back is either connected, or
splits into a pair of opposite components. We say that $C$ and the
corresponding hypersurfaces $X\in C$ are {\it reversible} in the
first case, and {\it irreversible} in the second one. In other
words, $X$ is reversible if its defining homogeneous polynomial,
$f$, can be continuously changed into $-f$ without creating
singularities in the process of deformation. One can extend the
notion of reversibility to singular varieties replacing
non-singular continuous families of equations by equisingular
families.

If the degree $d$ is even, then the region in $\Rp{n+1}$ where
$f>0$ defines a coorientation of $X_\R$, and reversibility
obviously means possibility to reverse this coorientation by a
deformation. If $n$ is odd, then such reversibility for
non-singular hypersurfaces is impossible, because the regions
where $f>0$ and $f<0$ are homologically different: they are
distinguished by the highest dimension in which the inclusion
homomorphism is nonzero. If $n$ is even, then reversibility is
possible: for example, a quadric is reversible if the signature of
its equation vanishes and irreversible otherwise. Furthermore, it
is not difficult to show that a real non-singular quartic surface
$X_\R$ is irreversible if it has more than one connected
components, as well as if it has a single component which is
contractible in $\Rp3$.
 Conversely, if $X_\R$ is connected and
non-contracible, then the quartic is reversible, at least if the
genus of $X_\R$ is $<10$ (the extremal case, $g=10$, remains
unknown to the authors). Thus we obtain nine reversible cases,
more than one hundred irreversible ones, and a unique unresolved
case.

If the degree $d$ is odd and $n$ is even, then $X_\R$ is
reversible for a trivial reason, because $-\id$ and $\id$ belong
to the same connected component of $GL(n+2,\R)$, and $f(-x)=-x$.
 If the both $d$ and $n$ are odd, then $f$ determines an
 orientation of $X_\R$ and reversibility obviously means
 possibility to alternate this orientation.
If $X_\R$ is symmetric with respect to a mirror reflection, then
such an alternation is realizable by a projective transformation,
which is one of manifestations of the similarity between the
notions of reversibility and achirality. Existence of symmetric
models proves in particular reversibility of curves of degree
$\le5$.

In the case of non-singular cubic threefolds the problem of
reversibility is already not trivial. The deformation
classification of such cubics obtained in \cite{Kr} gives 9
classes. Our analysis has shown that just one of these classes is
irreversible, namely, the class denoted
 $\Cal B(1)_I'$
in \cite{Kr}.


\Refs\widestnumber\key{ABCD}


\ref{ACT} \paper Real cubic surfaces and real hyperbolic geometry
\by D.~Allcock, J.~Carlson, and D. Toledo \jour CRAS \vol 337 \yr
2003 \pages 185--188
\endref\label{ACT}


\ref{B}
 \by N. Bourbaki
\book Groupes et Alg\`{e}bres de Lie, II
 \yr1968
 \bookinfo Hermann, Paris
 \pages
\endref\label{Bourbaki}

\ref{DIK} \paper Real Enriques surfaces \by A. Degtyarev,
I.Itenberg, V.Kharlamov \jour Lecture Notes Math., Springer \vol
1746 \yr 2000 \pages 259 pages
\endref\label{DIK}

\ref{D} \paper Reflection groups in algebraic geometry \by
I.~V.~Dolgachev \jour arXiv:math/0610938
\endref\label{Dolgachev}


\ref{FK}
 \by S. Finashin, V. Kharlamov
 \paper Deformation classes of real four-dimensional cubic
 hypersurfaces
 \jour J. Alg. Geom
 \vol
 \issue
 \yr
 \pages
\endref\label{FK}

\ref{Kh1}\by V.~Kharlamov \paper On classification of nonsingular
surfaces of degree 4 in ${\Bbb R} P^3$ with respect to rigid
isotopies \jour Funkt. Anal. i Priloz. \issue 1 \yr 1984 \pages
49--56
\endref\label{Kh1}

\ref{Kh2}\by V.~Kharlamov \paper On non-amphichaeral surfaces of
degree 4 in ${\Bbb R} P^3$ \jour Lecture Notes in Math. \vol 1346
\yr 1988 \pages 349--356
\endref\label{Kh2}

\ref{K} \paper \"{U}ber Flachen dritte Ordnung \by F.~Klein \jour
Math. Annalen \vol 6\yr 1873\pages 551--581
\endref\label{Kl}

\ref{K2} \paper \"{U}ber Flachen dritte Ordnung \by F.~Klein
\inbook Gesammelte Mathematische Abhandlungen \vol 2 \publ Verlag
von Julius Springer \yr 1921 \pages 11--67
\endref\label{Kl2}

\ref{Kr}
 \by V.~Krasnov
 \paper Rigid isotopy classification of real three-dimensional
 cubics
 \jour Izvestiya: Mathematics
 \vol 70
 \issue 4
 \yr2006
 \pages 731--768
\endref\label{Kr}



\ref{La}
 \by R. Laza
 \paper The moduli space of cubic fourfolds via the period map
 \jour arXiv:0705.0949
 \vol
 \issue
 \yr
 \pages
\endref\label{Laza}

\ref{Lo}
 \by E. Looijenga
 \paper The period map for cubic fourfolds
 \jour arXiv:0705.0951
 \vol
 \issue
 \yr
 \pages
\endref\label{Looijenga}


\ref{N1}
 \by V.~V.~Nikulin
 \paper Integer quadratic forms and some of their geometrical applications
 \jour Math. USSR -- Izv.
 \vol 43
  \yr 1979
 \pages 103--167
\endref\label{Nikulin}

\ref{N2}
 \by V.~V.~Nikulin
 \paper Remarks on connected components of moduli of real polarized K3 surfaces
 \jour arXiv:math.AG/0507197
 \yr 2006
\endref\label{N2}

\ref{R} \paper Zur Classification der Fl\"{a}chen dritter Ordnung
\by C.~Rodenberg \jour Math. Annalen \vol 14 \yr 1879\pages
46--110
\endref\label{R}


\ref{V}
 \by C.~Voisin
 \paper Th\'{e}or\`{e}me de Torelli pour les cubiques de $P^5$
 \jour Invent. Math.
 \vol 86
 \yr 1986
 \pages 577--601
\endref\label{Voisin}

\ref{Vin1}
 \by E. B. Vinberg
 \paper Some arithmetical discrete groups in Loba\v cevski\^{\i} spaces
 \inbook In: Proc. Int. Coll. on Discrete Subgroups of Lie Groups and Applications to Moduli
 (Bombay, 1973)
 \bookinfo Oxford University Press
 \yr 1975
 \pages 323--348
\endref\label{Vinberg1}



\ref{Vin2}
 \by E. B. Vinberg
 \paper Hyperbolic reflection groups
 \jour  Russian Math. Surveys
 \vol 40
 \yr 1985
 \pages 31--75
\endref\label{Vinberg2}


\ref{Vin3}
 \by E. B. Vinberg
 \paper The Two Most Algebraic K3 Surfaces
 \jour Math. Ann.
 \vol 265
 \yr 1983
 \pages 1--21
\endref\label{V3}

\endRefs
\enddocument